\documentclass{article}
\usepackage{amsfonts}
\usepackage{amsmath}
\usepackage{amssymb}
\usepackage{xcolor}

\newtheorem{theorem}{Theorem}

\newtheorem{corollary}[theorem]{Corollary}

\newtheorem{definition}[theorem]{Definition}

\newtheorem{lemma}[theorem]{Lemma}

\newtheorem{proposition}[theorem]{Proposition}
\newtheorem{remark}[theorem]{Remark}

\newenvironment{proof}[1][Proof]{\noindent\textbf{#1.} }{\ \rule{0.5em}{0.5em}}

\begin{document}

\title{Curve-Supported Concentrations of Vortices and Screw Dislocations Subject to Confinement}
\author{Carlos Garc\'{\i}a-Azpeitia\thanks{Departamento de Matem\'{a}ticas y Mec\'{a}nica, IIMAS-UNAM. Apdo. Postal 20-126, Col. San \'{A}ngel, Mexico City, 01000, Mexico. cgazpe@aries.iimas.unam.mx. ORCID: \mbox{0000-0002-6327-1444}.}, Thomas Hudson\thanks{Warwick Mathematics Institute, University of Warwick, Coventry, CV4 7AL, United Kingdom. ORCID: \mbox{0000-0002-8076-4937}. TH acknowledges partial support from UKRI during this work through NSF-UKRI grant reference
UKRI3611.}
}
\maketitle

\begin{abstract}

We study the concentration of logarithmically interacting defects, such as
Ginzburg--Landau vortices and screw dislocations, on closed curves in two dimensions, subject to a confining field. We prove that a measure supported on a
smooth closed curve solves the mean-field transport equation exactly when, in
Lagrangian labels, the curve moves with the negative induced velocity and each
label keeps its mass. This contour law is a gradient flow of the mean-field energy. We prove that the
gradient of the energy on a curve has a leading term proportional to the operator $|D|$ with Fourier symbol $|k|$, which prevents every smooth equilibrium curve from being a local minimizer, in a
gradient-flow analogue of the Kelvin--Helmholtz instability of vortex sheets. In a radially symmetric case with quadratic confinement,
we find an explicit circular equilibrium whose
linearization is diagonal in Fourier modes, and we show that, as $N\rightarrow \infty $, the linear spectrum of the regular $N$-gon converges, mode by mode, to that of the circle. Finally, we compute a critical point on two concentric circles
with lower energy than the single circle, indicating a
rich family of multi-circle equilibria.


\end{abstract}

\section{Introduction}

\label{sec:intro}

\subsection{Logarithmically interacting defects under confinement}

\label{subsec:intro-rings}

Large collections of interacting point defects (in this paper, as we discuss below, vortices in Ginzburg--Landau superconductors, screw dislocations, and vortices in fluids) can organize along closed
curves. In type-II superconductors, ring-shaped domains of vortices have been observed
experimentally~\cite{Groger2003}. Such a ring can be created by heating, and is metastable. In crystals, screw dislocations may accumulate in sheets called pile-ups, which then carry elastic incompatibility and can precede macroscopic failure, so that controlling such
accumulations is relevant to materials design~\cite{Kim19}. In fluids, rigidly rotating rings of point
vortices are a classical topic of study~\cite{Newton01}.
Motivated by these examples, we study such curve-supported concentrations within a single
mathematical model.

In all three settings described above a collection of defects can be modeled as points in two dimensions
interacting through a long-range logarithmic
potential, and a confining term keeps them together. In the appropriate asymptotic regime where
many well-spaced defects with the same topological charge are present and the confinement is sufficiently strong,
the empirical distribution of the defects converges to a measure $\mu $, with an internal
energy described by the functional
\begin{equation}
E[\mu ]=\frac{1}{2}\iint_{\mathbb{R}^{2}\times \mathbb{R}^{2}}G(x,y)\,d\mu
(x)\,d\mu (y)+\omega \int_{\mathbb{R}^{2}}V(x)\,d\mu (x),  \label{eq:E-intro}
\end{equation}%
where $G(x,y)=-\tfrac{1}{2\pi }\log |x-y|$ is the free-space Green's function
of $-\Delta $ and $V$ is a smooth confining potential of strength $\omega >0$: here, $\omega$ is a dimensionless constant that encodes the relative competition between the interaction term and the applied external potential.
The model example of a confining potential, and the one we use in all
explicit computations, is the quadratic trap
\begin{equation}
V(x)=\tfrac{1}{2}|x|^{2}.  \label{eq:V-quadratic}
\end{equation}%
The potential $V$ replaces boundary confinement, which in a bounded domain would be encoded in the
regular part of the Green's function. This quadratic trap induces
a drift $-\omega \nabla V=-\omega x$, that is, a restoring force pointing toward the origin, and growing linearly with distance.

The quadratic trap is a simple model of confinement appropriate for each of the three settings described. In the magnetic Ginzburg--Landau mean-field theory, such confinement can be generated by an applied field and the screening currents~\cite{SS2}. For screw dislocations,
a constant applied torsion produces the quadratic far-field confinement, and
$\omega $ is the torsion strength~\cite{Eshelby53}. For point
vortices, $\omega $ is the rotation frequency: the critical points
of~\eqref{eq:E-intro} are the configurations that rotate rigidly at angular
velocity $\omega $, the relative equilibria of the vortex problem.

In all three physical situations, curve-supported critical
points of~\eqref{eq:E-intro} are the central object of this
paper. For Ginzburg--Landau vortices and crystal dislocations the defects move
by the gradient flow of the energy, and it is this gradient evolution that we
develop. The dynamics of fluid vortices is instead Hamiltonian; we only
comment on the latter dynamics in Section~\ref{subsec:intro-contour}, and for fluid vortices
our results concern only steady configurations.

\subsection{Defect transport}

To account for the collective dynamics of point defects at finite temperature, we may adopt a
McKean--Vlasov (mean-field Langevin) perspective. At positive temperature, defects are subject
to Brownian diffusion, and the mean-field defect distribution is
absolutely continuous with respect to the Lebesgue measure, having density $\rho $ (so that $d\mu =\rho \,dx$). It
evolves through a competition between thermal diffusion and the
self-consistent nonlocal drift together with the external force, so that
\begin{equation}
  \label{eq:McKean} 
\begin{gathered}
\partial _{t}\rho =D\Delta \rho +\nabla \cdot \bigl(\rho \,v_{\rho }\bigr),
\\
\text{where}\qquad v_{\rho }(x,t) =\int_{\mathbb{R}^{2}}\nabla _{x}G(x,y)\,\rho
(y,t)\,dy+\omega \nabla V(x).
\end{gathered}
\end{equation}%
Formally, in the sense of~\cite{AGS08,CMV03}, equation~\eqref{eq:McKean} is
the $L^{2}$-Wasserstein gradient flow of the free-energy functional
$\mathcal{F}[\rho ]=D\int_{\mathbb{R}^{2}}\rho \log \rho \,dx+E[\rho ]$, with
$E$ as in~\eqref{eq:E-intro}, in which the entropy and the
interaction-plus-trap energy compete. For related aggregation--diffusion
equations, ground states with logarithmic interaction are analyzed
in~\cite{CCV15}, and long-time behavior and radial symmetry
in~\cite{CHVY19}.

While the diffusion term acts as an entropic regularizer that enforces
absolute continuity of the measure, the formal zero-temperature limit ($D\downarrow
0$) recovers a purely deterministic transport regime for the measure $\mu $,
\begin{equation}
\partial _{t}\mu =\nabla \cdot (\mu \,v_{\mu }).  \label{eq:transport}
\end{equation}%
We refer to this equation as the \emph{Eulerian} formulation of the defect dynamics of interest.
Unlike its diffusive counterpart, this transport law does not instantaneously smooth the defect distribution, so a solution initially supported on a curve may remain supported on a curve. This singular transport
regime is the primary focus of the present work. Without entropy, attractive
interaction kernels may even concentrate mass into singular measures in
finite time. The global existence of measure-valued solutions and the
occurrence of finite-time aggregation for nonlocal interaction equations with mildly singular interaction potentials (a class that excludes the logarithmic kernel) are established in~\cite{CDFLRS11}, and the related aggregation--diffusion equations, with their asymptotics and their vanishing-diffusion limit, are reviewed
in~\cite{CCY19}. A rigorous derivation of the deterministic mean-field limit
from the particle dynamics, for Coulomb and Riesz systems in
$\mathbb{R}^{d}$ without boundary effects, is provided in~\cite{Ser20}, under
the assumption that the limit solution is regular enough. That work focuses on the
boundary-free setting of the present paper, in which confinement is external
rather than geometric. Minimizers of a related anisotropic problem, among measures supported on ellipsoids, are studied in~\cite{Mora26}.

A related connection arises through fractional Riesz potentials. In the
singular limit $s\rightarrow 1$ of the inverse of $(-\Delta )^{s}$, one
recovers the transport equation~\eqref{eq:transport} with $\omega =0$,
governed by the Green's function of the Laplacian~\cite{CSV14}. In the
variational setting, \cite{BCLR13} relates the behavior of the interaction
kernel near the origin to the Hausdorff dimension of the support of local
minimizers of the interaction energy, and hence to when concentration on
lower-dimensional sets, such as rings, can occur.

The nonlocal transport equation we study belongs to a broad family of
continuity equations of the form~\eqref{eq:transport} arising in
distinct physical contexts. In superconductivity, the
mean-field model of Chapman, Rubinstein and Schatzman~\cite{CRS96} (see
also~\cite{E94}) transports the vortex density $\mu $ by the gradient of the
induced magnetic field $h$, which solves the London equation $-\Delta h+h=\mu
$ in the sample, with $h$ equal to the applied field, normalized to $1$, on its boundary. Its interaction kernel is
therefore the Green's function of the modified Helmholtz operator $(-\Delta
+\mathrm{Id})$, which has the logarithmic singularity of $G$ but decays
exponentially at large distances. Existence, and in some cases uniqueness, of solutions for this model and its
variants has been established in~\cite{DF18,MZ05,DZ03}.

In dislocation dynamics, screw dislocations interact through a kernel with the same logarithmic singularity. Their dynamics in the limit of many dislocations has been studied in~\cite{GvMPS20} (see~\cite{vMM19,vMPP22} for dislocations of both signs moving on a line, with annihilation), and a McKean--Vlasov equation like~\eqref{eq:McKean} is derived from an atomistic model in~\cite{Hudson20}.

\subsection{Curve-supported equilibria and instability}

\label{subsec:intro-contour}

One of our main contributions is the rigorous derivation of an evolution law for
curve-supported measures in a Lagrangian contour formulation, in contrast to the Eulerian form \eqref{eq:transport}. Let $\Gamma
_{t}$ be a time-dependent smooth closed curve parametrized by $\gamma (\sigma
,t)$ for $\sigma \in \mathbb{T}:=\mathbb{R}/2\pi \mathbb{Z}$, which carries a
line defect density $\lambda (\sigma ,t)$ (mass per unit parameter). The measure
$\mu_{t}=\lambda |\partial _{\sigma }\gamma |^{-1}\,\mathcal{H}^1\big| _{\Gamma_t}$ induces on the curve the velocity defined by
\begin{equation}
  \begin{aligned}
v_{\gamma ,\lambda }(\sigma ,t)&=\frac{1}{\lambda }\nabla _{\gamma }E[\mu
_{t}]\\
&=\mathrm{p.v.}\int_{\mathbb{T}}\nabla _{x}G\bigl(\gamma (\sigma
,t),\gamma (\sigma ^{\prime },t)\bigr)\,\lambda (\sigma ^{\prime
},t)\,d\sigma ^{\prime }+\omega \nabla V\bigl(\gamma (\sigma ,t)\bigr).
\end{aligned}
\label{eq:velocity}
\end{equation}
In Lagrangian coordinates each label $\sigma $ follows the flow and carries the
fixed mass $\lambda (\sigma )\,d\sigma $, so the contour law reads
\begin{equation}
\partial _{t}\gamma (\sigma ,t)=-v_{\gamma ,\lambda }(\sigma ,t),\qquad
\partial _{t}\lambda (\sigma ,t)=0.  \label{eq:lagrangian-intro}
\end{equation}%
The measures $\mu _{t}$ then
solve~\eqref{eq:transport}, and every smooth curve-supported solution of the equation can be
consistently relabeled into this form (Section~\ref{sec:lagrangian}). The total mass
$M=\int_{\mathbb{T}}\lambda \,d\sigma $ is constant, and the contour law is a
gradient flow of the energy (Section~\ref{subsec:gradient-flow}). Only the
normal component of $v_{\gamma ,\lambda }$ changes the shape of $\Gamma _{t}$,
while its tangential component simply redistributes mass along the curve. This
leads to the possible distinction between a \textit{geometric equilibrium}, in which normal velocities vanish, but tangential velocities may be nonzero,
and a \textit{Lagrangian equilibrium}, in which the total velocity on the curve vanishes. As we will see in the particular case of a circular equilibrium in a quadratic trap, a rigid rotation is a neutral mode of the
linearization; see Theorem~\ref{thm:circle-spectrum}.

We decompose the contour
velocity $v_{\gamma ,\lambda }$ into two components
(Section~\ref{sec:energy-curves}), these being (i) a density-coupled singular operator that recovers the
Fourier multiplier $\pi |D|$, where $|D|$ denotes the Fourier multiplier with
symbol $|k|$ acting on $2\pi $-periodic functions, and (ii) a remainder given by absolutely
convergent integrals, in which the leading-order short-range singularities
cancel. The external trap, being
smooth, enters only at order zero and does not affect the leading symbol.

Since the contour law is a gradient flow, its linearization at an equilibrium
is given by the second variation of the energy (Section~\ref{sec:geometric}),
and the leading term makes the second variation negative on
short-wavelength normal displacements, so that no smooth equilibrium curve is a local minimizer of the energy (Theorem~\ref{thm:energy-instability}). Every smooth equilibrium curve is
therefore formally unstable, with wrinkles whose amplitude grows
exponentially in time at a rate of order $|k|$ (Section~\ref{subsec:energy-instability}). This is the analogue, for the gradient-flow kernel, of the
Kelvin--Helmholtz instability of vortex sheets in fluid mechanics. This instability does not deprive the rings of physical relevance. As noted above in
Section~\ref{subsec:intro-rings}, such equilibria may persist for long enough to be observed experimentally, particularly if there are significant frictional effects which prevent their break-up.

On the other hand, this picture of instability is consistent with the variational intuition. The minimizer of
$E$ over all probability measures is the classical Frostman equilibrium
measure, which for the quadratic trap is a uniform density on a disk, of
radius $1$ at the matched strength $\omega _{\ast }=\frac{1}{2\pi}$ and unit
mass~\cite{ST97,SandierSerfaty2015}. Such curve-concentrated measures can thus at best approximate this equilibrium, for
instance through the generation of many concentric rings (Section~\ref{sec:numerics}). One therefore expects the energy to have
an infimum but no minimum among curve-supported measures, with minimizing
sequences that spread over the disk instead of staying on curves.

Our contour formulation is closely related to the Birkhoff--Rott equation
governing vortex sheets in fluid mechanics~\cite{Saffman92}, \cite[Section~9.2]{MB02}. The two
share the same integro-differential structure, but the fluid problem uses
the incompressible Biot--Savart kernel, which is the \emph{perpendicular}
gradient of the logarithmic potential, rather than the gradient kernel
appearing in~\eqref{eq:velocity}. Concretely, the fluid evolution is
obtained from ours by composing the velocity with the symplectic rotation
$J$ (a counterclockwise rotation of $\mathbb{R}^{2}$ through $\frac{\pi}{2}$),
\begin{equation}
v_{\gamma ,\lambda }^{\mathrm{fluid}}(\sigma ,t)=\frac{1}{\lambda }J\,\nabla
_{\gamma }E[\mu _{t}].\label{eq:velocity-fluid}
\end{equation}%
This turns the gradient field into a divergence-free one. The two laws therefore have opposite character. The vortex
sheet evolves by Hamiltonian dynamics, and its velocity produces
\emph{shear} across the curve, whereas the gradient kernel produces
\emph{stretching and compression}. Discretizing the sheet by point vortices
gives the classical Helmholtz--Kirchhoff dynamics of point
vortices~\cite{Newton01,Krasny86}, whose solutions rotating rigidly at
angular velocity $\omega $ are the critical points of the energy with the
quadratic trap of strength $\omega $. The
Birkhoff--Rott equation is one of the oldest and most intensively studied
problems in fluid mechanics, with a literature far too large to survey here.
As mentioned above, we do not focus on the fluid vortex case, and develop only the gradient
law~\eqref{eq:lagrangian-intro}. The contour-dynamics method for vortex
patches~\cite{ZHR79,Pullin92} uses the same Biot--Savart kernel to track
boundaries of regions of constant vorticity.

\subsection{Explicit equilibria and their discrete approximations}

To provide tractable examples of equilibria of the form we study, we analyze defects concentrated on a centered circle in a quadratic
trap. When the circle has uniform defect density, we show that it is an equilibrium of the curve evolution, and that its radius is the limit of the radii of the regular $N$-gon equilibria of the
discrete energy as $N\rightarrow \infty $. The circle is thus the continuum
analogue of the simplest relative equilibrium of the confined $N$-vortex
problem. The closest work in spirit to the present paper is
the dynamical analysis of closed curves in a Birkhoff--Rott framework
introduced in~\cite{SUB12}, which allows kernels with a gradient part, an incompressible part or both, and studies how the gradient part changes the Kelvin--Helmholtz instability.

On the circle, we show that the linearized flow is diagonal in Fourier modes
(Theorem~\ref{thm:circle-spectrum}). The normal displacement in the mode $k$
grows at the exact rate $\omega (|k|-2)$, so the modes with $|k|\geq 3$ grow
and so the circular equilibrium we find is unstable. On the other hand, the regular $N$-gon that approximates the circle behaves differently. It is stable for small $N$, because the
discreteness lowers the growth rates, most strongly those of the shortest waves (Section~\ref{subsec:linearization-circle}).

Finally, we compute the limit energy on measures supported on two concentric
circles. The optimal two-circle measure splits its mass between an inner ring
carrying roughly $40\%$ and an outer ring carrying roughly $60\%$, and has
strictly lower energy than the single circle. The curve-supported critical
points of the limit energy thus form a much richer family than the single
circle suggests, organized by the number of concentric rings and by their
mass partition.

\paragraph{Outline.}
The remainder of the paper is organized as follows: Section~\ref{sec:equilibria} presents a formal derivation of the mean-field energy from the discrete
model, and identifies the centered circle as the continuum limit
of the regular $N$-gon. Section~\ref{sec:energy-curves} characterizes the energy on curve-supported measures, computes its first variation, and isolates the leading singular operator in its gradient. Section~\ref{sec:lagrangian} proves that the
Eulerian transport law and the Lagrangian contour law are equivalent on
curve-supported measures, shows that Lagrangian labels can always be chosen,
and writes the contour law as a gradient flow. Section~\ref{sec:geometric} linearizes this
gradient flow at an equilibrium curve through the second variation of the
energy, computes the exact growth rates on the circle, compares them with those of the
regular polygon, and proves that no smooth equilibrium curve is a local minimizer of the energy. Section~\ref{sec:numerics} computes the limit energy on two
concentric circles, exhibits a two-ring critical point of lower energy, and
comments on configurations with more rings.

\section{Motivating examples}

\label{sec:equilibria}

In this section, we discuss some motivating examples of equilibria, showing that we can construct critical points of the energy $E$ which concentrate on curves; these serve as the motivation for our work on the rigorous definition of the energy which follows in subsequent sections.

\subsection{From the discrete energy to the mean-field energy}


For $N$ discrete defects at positions $x_{1},\dots ,x_{N}\in \mathbb{R}%
^{2} $, we can define a discrete energy
\begin{equation}
E_{N}(x_{1},\dots ,x_{N})=-\frac{1}{4\pi }\sum_{\substack{ j,k=1  \\ j\neq k
}}^{N}\log |x_{j}-x_{k}|+\omega N\sum_{j=1}^{N}V(x_{j}),  \label{eq:En}
\end{equation}%
which is consistent with the mean-field energy defined in \eqref{eq:E-intro}.
Here, in order that the confinement balance the interaction in the
mean-field limit, the trap is scaled by $N$ (equivalently, each particle
feels the mean trap generated by all $N$ of them).


The overdamped point defect dynamics given by the gradient descent of this energy with unit mobility is
\begin{equation}
\partial _{t}x_{j}=-\nabla _{x_{j}}E_{N}(x_{1},\dots ,x_{N}).
\label{eq:point-dynamics}
\end{equation}%
Differentiating~\eqref{eq:En} gives the explicit law
\begin{equation}
\partial _{t}x_{j}=\frac{1}{2\pi }\sum_{\substack{ k=1  \\ k\neq j}}^{N}%
\frac{x_{j}-x_{k}}{|x_{j}-x_{k}|^{2}}-\omega N\,\nabla V(x_{j}),
\label{eq:point-dynamics-expanded}
\end{equation}%
where Coulombic repulsion is balanced against the
external drift $-\omega N\,\nabla V$.

\paragraph{The continuum limit.}

Consider a sequence of configurations with $N$ defects at positions $x_{i}$%
, all with equal mass, and define the empirical measure
\begin{equation}
\mu _{N}=\frac{1}{N}\sum_{i=1}^{N}\delta _{x_{i}}.
\label{eq:empirical}
\end{equation}%
Assume that $\mu_{N}\rightharpoonup \mu$ as $N\to\infty$,
where the convergence is weak-$\ast$ convergence in the space of measures $\mathcal{M}(\mathbb{R}^2)$.
Dividing the energy $E_N$ by $N^{2}$ allows us to write both sums in the definition as
integrals against $\mu _{N}$. After subtracting the diagonal core energies,
the discrete energy converges to the mean-field Coulomb-gas functional
\begin{equation*}
E[\mu ]=\frac{1}{2}\iint_{\mathbb{R}^{2}\times \mathbb{R}^{2}}G(x,y)\,d\mu
(x)\,d\mu (y)+\omega \int_{\mathbb{R}^{2}}V(x)\,d\mu (x)
\end{equation*}
as defined in \eqref{eq:E-intro}. This convergence holds in the sense of $\Gamma $-convergence with respect to weak convergence of probability measures, and for the quadratic trap~\eqref{eq:V-quadratic} the empirical measures of minimizers of $E_{N}$ converge to the minimizer of $E$; see~\cite[Chapter~2]{Serfaty15}.

Physically, the functional~\eqref{eq:E-intro} balances logarithmic
repulsion against the trap $V$. The repulsion alone would spread the mass to
infinity, and the trap alone would collapse it to the origin. Their competition
selects an equilibrium measure of definite size. For the quadratic trap,
the unconstrained minimizer is the Frostman equilibrium measure, which is a uniform
density on a disk. As we show below, the energy also admits explicit
\emph{curve-supported} critical points (multi-circle configurations), which are the
equilibria of interest here.

\subsection{Circular equilibria as the limit of regular polygons}

\label{sec:circle-example}

We study the limit energy~\eqref{eq:E-intro} on its simplest
nontrivial critical point, the centered circle, and show that it is the $%
N\rightarrow \infty $ limit of the regular polygon minimizing the discrete
energy~\eqref{eq:En} among regular $N$-gons. Throughout this example we take the quadratic trap~%
\eqref{eq:V-quadratic}, with a general trap strength $\omega >0$.

\paragraph{Regular polygons.}
Place $N$ defects at the vertices of a regular $N$-gon of radius $r$,
\begin{equation*}
x_{j}=r\,e^{2\pi ij/N},\qquad j=0,\dots ,N-1.
\end{equation*}%
Two elementary identities reduce the energy~\eqref{eq:En} to a function of $%
r $ alone. The chordal distances are $|x_{i}-x_{j}|=2r\sin
\tfrac{\pi |i-j|}{N}$, and the product formula
\[
  \prod_{k=1}^{N-1}\bigl(2\sin \tfrac{%
    \pi k}{N}\bigr)=N\]
allows us to write
\begin{equation*}
\sum_{\substack{ i,j=1  \\ i\neq j}}^N\log |x_{i}-x_{j}|=N\sum_{k=1}^{N-1}\log %
\Bigl(2r\sin \tfrac{\pi k}{N}\Bigr)=N(N-1)\log r+N\log N.
\end{equation*}%
Secondly, we note that $\sum_{j=1}^NV(x_{j})=\tfrac{1}{2}\sum_{j}|x_{j}|^{2}=\tfrac{1}{2}Nr^{2}$%
. Substituting these expressions into~\eqref{eq:En}, we can express the energy as a function of the radius $r$ alone:
\begin{equation}
E_{N}(r)=-\frac{1}{4\pi }\bigl[N(N-1)\log r+N\log N\bigr]+\frac{\omega N^{2}%
}{2}r^{2}.  \label{eq:EN-polygon}
\end{equation}%
Differentiating and discarding the $r$-independent term gives $E_{N}^{\prime
}(r)=-\tfrac{N(N-1)}{4\pi r}+\omega N^{2}r$, so the unique stationary radius
in $(0,\infty )$ is
\begin{equation}
r_{N}=\sqrt{\frac{N-1}{4\pi \omega N}}.  \label{eq:rN}
\end{equation}%
Since $E_{N}^{\prime \prime }(r)>0$, this $r_{N}$ is the unique minimizer of
the energy restricted to regular $N$-gons.
Finally, we note that this radius does indeed induce an equilibrium configuration: by reflection symmetry, all angular components of the forces acting on the defects must vanish, and $r_N$ has been chosen to fix the radial force components to be zero.

\paragraph{The circle.} Let us now consider the continuum limiting case. Let $\nu _{r}$ denote the uniform probability measure
on the centered circle of radius $r$, and take $\mu =M\nu _{r}$ with total
mass $M$. Using this definition, we can compute
\begin{equation}
\iint \log |x-y|\,d\nu _{r}(x)\,d\nu _{r}(y)=\frac{1}{2\pi }\int_{0}^{2\pi
}\log \bigl|r\,e^{i\alpha }-r\bigr|\,d\alpha =\log r,  \label{eq:circle-self}
\end{equation}%
where we use the fact that $\tfrac{1}{2\pi }\int_{0}^{2\pi }\log |e^{i\alpha }-1|\,d\alpha =0$, which follows from the mean value property of the harmonic function $z\mapsto \log |z-1|$ on the circles $|z|=q<1$ and dominated convergence as $q\rightarrow 1$.
Since $\int V\,d\nu _{r}=\tfrac{1}{2}\int |x|^{2}\,d\nu _{r}=\tfrac{1}{2}r^{2}
$, the limit energy~\eqref{eq:E-intro} of $\mu =M\nu _{r}$ is again expressible as a function of $r$ only,
\begin{equation}
E(r)=-\frac{M^{2}}{4\pi }\log r+\frac{\omega M}{2}r^{2}.
\label{eq:Einf-circle}
\end{equation}%
The derivative $E^{\prime }(r)=-\tfrac{M^{2}}{4\pi r}+\omega M\,r%
$ vanishes at
\begin{equation}
r_{\ast }=\sqrt{\frac{M}{4\pi \omega }},\label{eq:rstar}
\end{equation}%
a strict local minimum among these competitors since $E^{\prime \prime }(r_{\ast })>0$. Note that the
equilibrium radius grows with the total mass $M$, as the repulsive interactions scale as $M^{2}$ and
the effect of the trap only as $M$, while increasing $\omega$ reduces the radius by pushing the defects inward. For comparison with the polygonal case above, it is natural to take the unit-mass
normalization $M=1$, since the empirical measure $\mu _{N}$ has mass one, and in this case $r_{\ast }$ is
exactly the $N\rightarrow \infty $ limit~\eqref{eq:rN} of the polygon
radius, $r_{N}$. Thus the regular $N$-gon minimizing the discrete energy among regular $N$-gons
converges to the uniform circle of radius $r_{\ast }$ minimizing the limit
energy among centered circles. The circle is also a critical point of $E$ for every variation of the curve at fixed density, because its velocity vanishes identically, as we show rigorously in Theorem~\ref{thm:circle-spectrum}.

\section{Making sense of the energy on curves}
\label{sec:energy-curves}
Motivated by the examples in the previous section, we now characterize the energy \eqref{eq:E-intro} as defined on curve-supported measures.

\subsection{Energy and first variation of a curve-supported measure}

\label{subsec:first-variation}

Let $\Gamma \subset \mathbb{R}^{2}$ be a smooth closed curve parametrized by
a smooth embedding $\gamma \colon \mathbb{T}\rightarrow \mathbb{R}^{2}$,
where $\mathbb{T}=\mathbb{R}/2\pi \mathbb{Z}$. Let $\lambda \colon \mathbb{T}%
\rightarrow \lbrack 0,\infty )$ be a smooth, nonnegative line density (mass
per unit parameter). The measure concentrated on $\Gamma $ with density $%
\lambda $ is
\begin{equation}
  \mu =\lambda|\gamma ^{\prime }|^{-1}\,\mathcal{H}^1\big|_{\Gamma},  \label{eq:curveAnsatzStatic}
\end{equation}%
where $\mathcal{H}^1$ is the one-dimensional Hausdorff measure, so that for every continuous function $\phi $ we have $\int \phi \,d\mu =\int_{\mathbb{T}}\phi (\gamma )\,\lambda \,d\sigma $.
The measure $\mu $ in~\eqref{eq:curveAnsatzStatic} does not depend on how we may choose to parametrize $\Gamma $; in fact, assuming sufficient regularity of the curve, it is a normal current.

A natural choice for parametrizing the measure is to take a parametrization proportional to arc length. Writing $\gamma$ and $\lambda$ as functions of this parameter $\sigma$, we have
\begin{equation}
\begin{split}
E[\mu]& =\frac{1}{2}\iint_{\mathbb{T}\times \mathbb{T}}G\!\bigl(%
\gamma (\sigma ),\gamma (\sigma ^{\prime })\bigr)\,\lambda (\sigma
)\,\lambda (\sigma ^{\prime })\,d\sigma \,d\sigma ^{\prime } +\omega \int_{\mathbb{T}}V\bigl(\gamma (\sigma )\bigr)\,\lambda (\sigma )\,d\sigma.
\end{split}
\label{eq:curve-energy-full}
\end{equation}%
It is clear that any alternative regular parametrization of the curve will yield an identical result, simply through an application of the change of variables formula. We write $E_{\lambda }(\gamma )$ for the right-hand side of~\eqref{eq:curve-energy-full}, regarded as a function of the curve $\gamma $ at fixed $\lambda $.
If we vary $\gamma \mapsto \gamma + \varepsilon h$ with a smooth perturbation $h \colon \mathbb{T}
\to \mathbb{R}^2$, the first variation of~%
\eqref{eq:curve-energy-full} is
\begin{equation}  \label{eq:varE-final}
\delta E[\mu][h] = \int_{\mathbb{T}} \lambda(\sigma) \Biggl( \int_{%
\mathbb{T}} \nabla_{x} G\bigl(\gamma(\sigma), \gamma(\sigma^{\prime })\bigr)%
\, \lambda(\sigma^{\prime })\, d\sigma^{\prime } + \omega \nabla V\bigl(%
\gamma(\sigma)\bigr) \Biggr) \cdot h(\sigma)\, d\sigma,
\end{equation}
where we invoke the symmetry $\nabla_{y} G(x,y) =
\nabla_{x} G(y,x)$. Consequently, defining the unweighted $L^{2}(\mathbb{T})$ gradient $\nabla_{\gamma}
E_{\lambda}(\gamma)$ through the duality $\delta E_{\lambda}[h] = \int_{%
\mathbb{T}} \nabla_{\gamma} E_{\lambda}(\gamma)(\sigma) \cdot h(\sigma)\,
d\sigma$, we obtain
\begin{equation}  \label{eq:gradX-final}
\nabla_{\gamma} E_{\lambda}(\gamma)(\sigma) = \lambda(\sigma) \Biggl( \int_{%
\mathbb{T}} \nabla_{x} G\bigl(\gamma(\sigma), \gamma(\sigma^{\prime })\bigr)%
\, \lambda(\sigma^{\prime })\, d\sigma^{\prime } + \omega \nabla V\bigl(%
\gamma(\sigma)\bigr) \Biggr).
\end{equation}
Using the definition of $G$ and interpreting singular integrals in the
principal-value sense, equation \eqref{eq:gradX-final} becomes
\begin{equation}
\nabla _{\gamma }E_{\lambda }(\gamma )(\sigma )=\lambda (\sigma )\Biggl(-%
\frac{1}{2\pi }\,\mathrm{p.v.}\!\int_{\mathbb{T}}\frac{\gamma (\sigma
)-\gamma (\sigma ^{\prime })}{\lvert \gamma (\sigma )-\gamma (\sigma
^{\prime })\rvert ^{2}}\,\lambda (\sigma ^{\prime })\,d\sigma ^{\prime
}+\omega \nabla V\bigl(\gamma (\sigma )\bigr)\Biggr).
\label{eq:gradX-decomp-final}
\end{equation}
These computations are legitimate: for a smooth embedding $\gamma $, the ratio $|\gamma (\sigma )-\gamma (\sigma ^{\prime })|^{2}/(4\sin ^{2}\frac{\sigma -\sigma ^{\prime }}{2})$ extends to a smooth positive function on $\mathbb{T}\times \mathbb{T}$, equal to $|\partial _{\sigma }\gamma (\sigma )|^{2}$ on the diagonal, so $\log |\gamma (\sigma )-\gamma (\sigma ^{\prime })|$ is $\log |2\sin \frac{\sigma -\sigma ^{\prime }}{2}|$ plus half the logarithm of this ratio, and the first term does not depend on the curve. In the same way the kernel in~\eqref{eq:gradX-decomp-final} is a smooth multiple of $\cot \frac{\sigma -\sigma ^{\prime }}{2}$ plus a smooth kernel, so the energy and its gradient depend smoothly on smooth perturbations of the curve and of the density.


\subsection{The leading singular operator of the gradient}

We now work with an arbitrary regular parametrization $\gamma \colon \mathbb{T}\rightarrow \mathbb{R}^{2}$ of the curve $\Gamma $, so that $|\gamma ^{\prime }(\sigma )|>0$ for every $\sigma $, and we write $\rho =\lambda /|\gamma ^{\prime }|$ for the density per unit arc length.

To isolate the principal one-dimensional singular operator, we introduce the
following three operators, with singular integrals understood in the
principal-value sense:
\begin{align}
\mathcal{S}_{\lambda }(\gamma )(\sigma )& :=\frac{\lambda (\sigma )}{|\gamma ^{\prime }(\sigma )|^{2}}\,%
\mathrm{p.v.}\!\int_{\mathbb{T}}\frac{\gamma (\sigma )-\gamma (\sigma
^{\prime })}{4\sin ^{2}\!\tfrac{\sigma -\sigma ^{\prime }}{2}}\,d\sigma
^{\prime },  \label{eq:S-def} \\
\mathcal{R}_{\lambda }(\gamma )(\sigma )& :=\mathrm{p.v.}\!\int_{\mathbb{T}}%
\frac{\gamma (\sigma )-\gamma (\sigma ^{\prime })}{\lvert \gamma (\sigma
)-\gamma (\sigma ^{\prime })\rvert ^{2}}\,\lambda (\sigma ^{\prime
})\,d\sigma ^{\prime }-\mathcal{S}_{\lambda }(\gamma )(\sigma ),
\label{eq:R-def} \\
\mathcal{V}(\gamma )(\sigma )& :=\omega \,\nabla V\bigl(\gamma (\sigma )%
\bigr).\vphantom{\int_{\mathbb{T}}} \label{eq:K-def}
\end{align}%
The operator $\mathcal{S}_{\lambda }$ freezes $\lambda $ and the speed $|\gamma ^{\prime }|$ at $\sigma $ and
uses the exact periodic kernel $4\sin ^{2}\!\tfrac{\sigma -\sigma ^{\prime }%
}{2}$, the squared chordal distance on $\mathbb{T}$. The operator $\mathcal{R%
}_{\lambda }$ is the remainder obtained by subtracting $\mathcal{S}_{\lambda
}$ from the full singular integral, and absorbs both the geometric mismatch
between $|\gamma (\sigma )-\gamma (\sigma ^{\prime })|^{2}$ and $|\gamma ^{\prime }(\sigma )|^{2}\cdot 4\sin ^{2}\!\tfrac{\sigma -\sigma ^{\prime }}{2}$, and the variation of $%
\lambda $ along $\Gamma $. The term $\mathcal{V}$ is the smooth external
drift. It is a local multiplication operator given by the gradient of the trap
along the curve, with no integral kernel at all. With these definitions, equation~\eqref{eq:gradX-decomp-final} can be rewritten as
\begin{equation}
\nabla _{\gamma }E_{\lambda }(\gamma )(\sigma )=\lambda (\sigma )\Bigl[-%
\frac{1}{2\pi }\bigl(\mathcal{S}_{\lambda }(\gamma )(\sigma )+\mathcal{R}%
_{\lambda }(\gamma )(\sigma )\bigr)+\mathcal{V}(\gamma )(\sigma )\Bigr].
\label{eq:grad-decomp-final-noarray}
\end{equation}

For convenience throughout the rest of our analysis, we make the following notational definition:

\begin{definition}[Differences]
\label{def:delta}We write $\delta \sigma =\sigma -\sigma ^{\prime }$ and, for a function $f$ on $\mathbb{T}$, $\delta f=f(\sigma )-f(\sigma ^{\prime })$. In particular, $\delta \gamma =\gamma (\sigma )-\gamma (\sigma ^{\prime })$ is the chord of the curve and $\delta \lambda =\lambda (\sigma )-\lambda (\sigma ^{\prime })$. A function written without its argument is evaluated at $\sigma $.
\end{definition}

\noindent
The following result now provides a Fourier characterization of the first operator, $\mathcal{S}_\lambda$.

\begin{lemma}
\label{lem:piD}For every smooth function $f$ on $\mathbb{T}$, with Fourier coefficients $\widehat{f}_{k}$,
\begin{equation}
\mathrm{p.v.}\!\int_{\mathbb{T}}\frac{f(\sigma )-f(\sigma ^{\prime })}{4\sin
^{2}\!\tfrac{\sigma -\sigma ^{\prime }}{2}}\,d\sigma ^{\prime }=\pi \lvert
D\rvert f(\sigma )=\sum_{k\in \mathbb{Z}}\pi \lvert k\rvert \,\widehat{f}%
_{k}\,e^{ik\sigma }.  \label{eq:piD-identity}
\end{equation}%
\end{lemma}

\begin{proof}
First let $f(\sigma )=e^{ik\sigma }$ with $k\in \mathbb{Z}$. Then
\begin{equation*}
\delta f=e^{ik\sigma }\bigl(1-e^{-ik\,\delta \sigma }\bigr)=e^{ik\sigma }\bigl(1-\cos k\,\delta \sigma \bigr)+i\,e^{ik\sigma }\sin k\,\delta \sigma .
\end{equation*}%
By definition, the principal value removes the interval $|\delta \sigma |<\varepsilon $ and lets $\varepsilon \rightarrow 0$. The second term is odd in $\delta \sigma $, so its integral over $\varepsilon <|\delta \sigma |\leq \pi $ vanishes for every $\varepsilon $. In the first term we use $1-\cos k\,\delta \sigma =2\sin ^{2}\frac{k\,\delta \sigma }{2}$, which gives
\begin{equation*}
\frac{1-\cos k\,\delta \sigma }{4\sin ^{2}\frac{\delta \sigma }{2}}=\frac{1}{2}\,\frac{\sin ^{2}\frac{k\,\delta \sigma }{2}}{\sin ^{2}\frac{\delta \sigma }{2}}.
\end{equation*}%
This function is nonnegative and bounded by $k^{2}/2$, so the principal value of its integral is an ordinary integral, and for $k=0$ it vanishes. For $k\neq 0$, the geometric sum $\sum_{j=0}^{|k|-1}e^{ij\,\delta \sigma }=(e^{i|k|\,\delta \sigma }-1)/(e^{i\,\delta \sigma }-1)$ has modulus $|\sin \frac{|k|\,\delta \sigma }{2}|/|\sin \frac{\delta \sigma }{2}|$, because $|e^{ix}-1|=2|\sin \frac{x}{2}|$. Since the sum has $|k|$ terms with coefficient $1$, Parseval's identity gives
\begin{equation}
\int_{\mathbb{T}}\frac{\sin ^{2}\frac{k\,\delta \sigma }{2}}{\sin ^{2}\frac{\delta \sigma }{2}}\,d\sigma ^{\prime }=\int_{\mathbb{T}}\Bigl|\sum_{j=0}^{|k|-1}e^{ij\,\delta \sigma }\Bigr|^{2}d\sigma ^{\prime }=2\pi |k|,  \label{eq:fejer}
\end{equation}%
so the left side of~\eqref{eq:piD-identity} equals $\frac{1}{2}\cdot 2\pi |k|\,e^{ik\sigma }=\pi |k|\,e^{ik\sigma }$, which is the right side.

More generally, for smooth $f$ we have $\sum_{k}|k|\,|\widehat{f}_{k}|<\infty $. On the set $|\delta \sigma |>\varepsilon $ the Fourier series of the integrand converges uniformly, so we may integrate term by term. By the argument above, the integral of the mode $k$ over this set is $e^{ik\sigma }$ times the integral of the nonnegative function above, so its modulus is at most $\pi |k|$. Dominated convergence for series then allows us to let $\varepsilon \rightarrow 0$ in each term, which gives~\eqref{eq:piD-identity}.
\end{proof}
\medskip

By Lemma~\ref{lem:piD} applied to the components of $\gamma $, the operator $\mathcal{S}_{\lambda }$ reduces \emph{exactly} to a
multiplication by $\lambda (\sigma )/|\gamma ^{\prime }(\sigma )|^{2}$ composed with the Fourier multiplier $\pi |D|$. So for every parametrization with $|\gamma ^{\prime }(\sigma )|>0$,
\begin{equation}
\mathcal{S}_{\lambda }(\gamma )(\sigma )=\frac{\pi \lambda (\sigma )}{|\gamma ^{\prime }(\sigma )|^{2}}%
\,\lvert D\rvert \gamma (\sigma ).  \label{eq:S-exact}
\end{equation}%
The following proposition shows that the integral defining $\mathcal{R}_{\lambda }(\gamma )$ needs no principal value, so that $\mathcal{S}_{\lambda }$ carries the whole singularity of the integral in~\eqref{eq:R-def}.

\begin{proposition}
\label{prop:R-regular} Let $\gamma \in C^{2}(\mathbb{T};\mathbb{R}^{2})$ be
a regular parametrization as before with $|\gamma ^{\prime }|>0$   and $%
\lambda \in C^{1}(\mathbb{T})$. Then the integrand of $\mathcal{R}_{\lambda }(\gamma )$ in~\eqref{eq:R-def} extends to a continuous function on $\mathbb{T}\times \mathbb{T}$. Hence $\mathcal{R}_{\lambda }(\gamma )$ is given by an absolutely convergent integral and is a continuous function on $\mathbb{T}$, and the whole principal-value singularity of~\eqref{eq:grad-decomp-final-noarray} is carried by $\mathcal{S}_{\lambda }(\gamma )$.
\end{proposition}

\begin{proof}
Write $\lambda (\sigma ^{\prime })=\lambda (\sigma )+[\lambda (\sigma
^{\prime })-\lambda (\sigma )]$ inside the full singular integral of~%
\eqref{eq:R-def} and, using the definition of $\mathcal{S}_{\lambda }$,
split $\mathcal{R}_{\lambda }(\gamma )=\mathcal{R}_{\lambda }^{(1)}(\gamma )+%
\mathcal{R}_{\lambda }^{(2)}(\gamma )$, where
\begin{align}
\mathcal{R}_{\lambda }^{(1)}(\gamma )(\sigma )& :=\lambda (\sigma )\,\mathrm{%
p.v.}\!\int_{\mathbb{T}}\bigl(\gamma (\sigma )-\gamma (\sigma ^{\prime })%
\bigr)  \notag \\
& \qquad \times \left[\frac{1}{\lvert \gamma (\sigma )-\gamma (\sigma ^{\prime
})\rvert ^{2}}-\frac{1}{|\gamma ^{\prime }(\sigma )|^{2}\cdot 4\sin ^{2}\!\tfrac{\sigma -\sigma
^{\prime }}{2}}\right]d\sigma ^{\prime },  \label{eq:R1-def} \\
\mathcal{R}_{\lambda }^{(2)}(\gamma )(\sigma )& :=\mathrm{p.v.}\!\int_{%
\mathbb{T}}\frac{\gamma (\sigma )-\gamma (\sigma ^{\prime })}{\lvert \gamma
(\sigma )-\gamma (\sigma ^{\prime })\rvert ^{2}}\,\bigl[\lambda (\sigma
^{\prime })-\lambda (\sigma )\bigr]\,d\sigma ^{\prime }.  \label{eq:R2-def}
\end{align}%
Both integrands are continuous away from the diagonal $\sigma ^{\prime
}=\sigma$, and we compute their limits as $\sigma ^{\prime }\to \sigma $.

\smallskip \noindent \textit{Kernel of $\mathcal{R}_{\lambda }^{(1)}$.}
A Taylor expansion of $\gamma $ at $\sigma $ gives $\delta \gamma =\gamma ^{\prime }\,\delta \sigma -\frac{1}{2}\gamma ^{\prime \prime }(\delta \sigma )^{2}+o((\delta \sigma )^{2})$. Since $|\gamma ^{\prime }|>0$, this entails
\begin{align*}
|\delta \gamma |^{2}& =|\gamma ^{\prime }|^{2}(\delta \sigma )^{2}-(\gamma ^{\prime }\cdot \gamma ^{\prime \prime })\,(\delta \sigma )^{3}+o((\delta \sigma )^{3}), \\
\frac{1}{|\delta \gamma |^{2}}& =\frac{1}{|\gamma ^{\prime }|^{2}(\delta \sigma )^{2}}+\frac{\gamma ^{\prime }\cdot \gamma ^{\prime \prime }}{|\gamma ^{\prime }|^{4}\delta \sigma }+o(1/\delta \sigma )
\end{align*}%
as $\delta \sigma \to0$.
Since $1/(4\sin ^{2}\frac{\delta \sigma }{2})=1/(\delta \sigma )^{2}+O(1)$, the bracket in~\eqref{eq:R1-def} is
\[
B:= \left[\frac{1}{\lvert \delta \gamma \rvert ^{2}}-\frac{1}{|\gamma ^{\prime }(\sigma )|^{2}\cdot 4\sin ^{2}\!\tfrac{\delta \sigma }{2}}\right]=\frac{\gamma ^{\prime }\cdot \gamma ^{\prime \prime }}{|\gamma ^{\prime }|^{4}\delta \sigma }+o(1/\delta \sigma ),
\]
where we note that the terms of order $1/(\delta \sigma )^{2}$ cancel because $\mathcal{S}_{\lambda }$ uses the local speed $|\gamma ^{\prime }(\sigma )|$. Multiplying by $\delta \gamma =\gamma ^{\prime }\,\delta \sigma +O((\delta \sigma )^{2})$, we find that the kernel of $\mathcal{R}_{\lambda }^{(1)}$ satisfies
\begin{equation*}
\lim_{\sigma ^{\prime }\rightarrow \sigma }\lambda (\sigma )\,\delta \gamma \,B=\frac{\lambda (\sigma )\,\bigl(\gamma ^{\prime }(\sigma )\cdot \gamma ^{\prime \prime }(\sigma )\bigr)}{|\gamma ^{\prime }(\sigma )|^{4}}\,\gamma ^{\prime }(\sigma ).
\end{equation*}%

\smallskip \noindent \textit{Kernel of $\mathcal{R}_{\lambda }^{(2)}$.}
Taylor expanding $\gamma $ and $\lambda $ at $\sigma $ gives
\begin{equation*}
-\frac{\delta \gamma \,\delta \lambda }{|\delta \gamma |^{2}}\;\xrightarrow[\sigma ^{\prime }\to \sigma ]{}\;-%
\frac{\gamma ^{\prime }(\sigma )\,\lambda ^{\prime }(\sigma )}{|\gamma ^{\prime }(\sigma )|^{2}}.
\end{equation*}

\smallskip\noindent
\textit{Conclusion.} These Taylor expansions hold uniformly in $\sigma $, because $\gamma ^{\prime \prime }$ and $\lambda ^{\prime }$ are uniformly continuous on $\mathbb{T}$. Hence both kernels extend continuously across the diagonal,
so no principal value is needed. Since they are continuous on the compact set $\mathbb{T}\times \mathbb{T}$, the integrals in~\eqref{eq:R1-def} and~\eqref{eq:R2-def} converge absolutely and define continuous functions of $\sigma $. The integrand of $\mathcal{R}_{\lambda }(\gamma )$ is the sum of these two integrands, which proves the proposition.
\end{proof}\medskip

\noindent Substituting~\eqref{eq:S-exact} into~%
\eqref{eq:grad-decomp-final-noarray} yields the gradient in its definitive
parametric form:
\begin{equation}
\nabla _{\gamma }E_{\lambda }(\gamma )(\sigma )=-\frac{\lambda (\sigma )^{2}%
}{2|\gamma ^{\prime }(\sigma )|^{2}}\,\lvert D\rvert \gamma (\sigma )+\lambda (\sigma )\Bigl[-\frac{1}{%
2\pi }\mathcal{R}_{\lambda }(\gamma )(\sigma )+\mathcal{V}(\gamma )(\sigma )%
\Bigr].  \label{eq:grad-final-general}
\end{equation}%
In terms of the density per unit arc length $\rho =\lambda /|\gamma ^{\prime }|$ introduced above, we can write the leading term as $-\frac{1}{2}\rho (\sigma )^{2}\lvert D\rvert \gamma $. 

Once the dynamics on curves has been established in Section~\ref{sec:lagrangian}, we will see in Section~\ref{sec:geometric} that the leading term of this decomposition is responsible for the instability of equilibrium curves.

\section{The contour law that governs curve-supported measures}

\label{sec:lagrangian}

In this section, we make precise the relation between the Eulerian transport law~\eqref{eq:transport} and the Lagrangian contour law~\eqref{eq:lagrangian-intro}, for measures carried by a curve: we show that the two laws are equivalent, that Lagrangian labels can always be chosen, and that in these labels the contour law is a gradient flow of the energy.

\subsection{The Lagrangian contour law}

\label{subsec:contour-law}

Assume that $\mu _{t}$ is supported on a smooth, time-dependent curve $%
\Gamma _{t}$, represented by a line density $\lambda (\sigma ,t)$ as in~\eqref{eq:curveAnsatzStatic}, that is,
\begin{equation}
\langle \mu _{t},\phi \rangle =\int_{\mathbb{T}}\lambda (\sigma ,t)\,\phi \bigl(\gamma (\sigma ,t)\bigr)\,d\sigma \qquad \text{for every }\phi \in C_{c}(\mathbb{R}^{2}).  \label{eq:curveAnsatz}
\end{equation}
Relabeling the curve by an orientation-preserving diffeomorphism $\tilde{\sigma}=\psi (\sigma )$ of $\mathbb{T}$ changes $\lambda $, but not $\mu _{t}$, provided that the mass of every arc is preserved, $\tilde{\lambda}(\tilde{\sigma})\,d\tilde{\sigma}=\lambda (\sigma )\,d\sigma $; the density per unit arc length $\rho =\lambda /|\partial _{\sigma }\gamma |$, on the other hand, is intrinsic.

The measure $\mu _{t}$ induces on $\Gamma _{t}$ the velocity $v_{\gamma ,\lambda }$ given by~\eqref{eq:velocity}, with the self-interaction understood in the principal-value sense. In Lagrangian labels each parameter value $\sigma $ is attached to a material particle, which carries the fixed mass $\lambda (\sigma )\,d\sigma $ and moves with velocity $-v_{\gamma ,\lambda }(\sigma ,t)$; this is the contour law~\eqref{eq:lagrangian-intro}.

\subsection{The contour law is equivalent to the transport law}

\label{subsec:equivalence}

A family of finite measures $(\mu _{t})_{t\in \lbrack 0,T)}$ solves~%
\eqref{eq:transport} weakly if, for every $\phi \in C_{c}^{\infty }(%
\mathbb{R}^{2})$, the map $t\mapsto \langle \mu _{t},\phi \rangle $ is
differentiable and
\begin{equation}
\frac{d}{dt}\langle \mu _{t},\phi \rangle =-\int_{\mathbb{R}^{2}}\nabla \phi
(x)\cdot v_{\mu }(x,t)\,d\mu _{t}(x).  \label{eq:weak-transport}
\end{equation}%
This identity follows from~\eqref{eq:transport} by one integration by
parts. For a measure carried by a curve, $v_{\mu }$ is evaluated on the curve
by the principal-value formula~\eqref{eq:velocity}. We identify $\mathbb{R}^{2}$ with $\mathbb{C}$ and orient each closed curve $\gamma (\cdot ,t)$ counterclockwise. We write $\tau
=\partial _{\sigma }\gamma /|\partial _{\sigma }\gamma |$ for the unit
tangent, $n=-i\tau $ for the outward unit normal, and
\begin{equation}
g(\sigma ,t):=\bigl(\partial _{t}\gamma (\sigma ,t)+v_{\gamma ,\lambda
}(\sigma ,t)\bigr)\cdot \tau (\sigma ,t)  \label{eq:gauge-function}
\end{equation}%
for the tangential velocity of the labels relative to the flow.
\begin{proposition}[Curve-supported weak solutions]
\label{prop:euler-lagrange}Let $T>0$. Let $\gamma \colon \mathbb{T}\times
\lbrack 0,T)\rightarrow \mathbb{R}^{2}$ and $\lambda \colon \mathbb{T}\times
\lbrack 0,T)\rightarrow (0,\infty )$ be smooth, with $\gamma (\cdot ,t)$ an
embedding for every $t$. Assume that the velocity $v_{\gamma ,\lambda }$
of~\eqref{eq:velocity} is smooth on $\mathbb{T}\times \lbrack 0,T)$,
and let $\mu _{t}$ be given by~\eqref{eq:curveAnsatz}. Then $\mu _{t}$
solves~\eqref{eq:transport} weakly if and only if
\begin{equation}
\partial _{t}\gamma \cdot n=-v_{\gamma ,\lambda }\cdot n  \label{eq:normal-law}
\end{equation}%
and
\begin{equation}
\partial _{t}\lambda =\partial _{\sigma }\Bigl(\frac{\lambda \,g}{|\partial
_{\sigma }\gamma |}\Bigr)  \label{eq:continuity-lambda}
\end{equation}%
hold on $\mathbb{T}\times (0,T)$.
\end{proposition}
\begin{proof}
Fix $\phi \in C_{c}^{\infty }(\mathbb{R}^{2})$ and $t$. The integrand in~\eqref{eq:curveAnsatz} is smooth and $\mathbb{T}
$ is compact, so we may differentiate under the integral sign. The chain rule
gives
\begin{equation*}
\frac{d}{dt}\langle \mu _{t},\phi \rangle =\int_{\mathbb{T}}\partial
_{t}\lambda \,\phi (\gamma )\,d\sigma +\int_{\mathbb{T}}\lambda \,\nabla \phi
(\gamma )\cdot \partial _{t}\gamma \,d\sigma .
\end{equation*}%
The right-hand side of~\eqref{eq:weak-transport} equals $-\int_{\mathbb{T}%
}\lambda \,\nabla \phi (\gamma )\cdot v_{\gamma ,\lambda }\,d\sigma $. Hence~%
\eqref{eq:weak-transport} holds if and only if
\begin{equation}
\int_{\mathbb{T}}\partial _{t}\lambda \,\phi (\gamma )\,d\sigma +\int_{%
\mathbb{T}}\lambda \,\nabla \phi (\gamma )\cdot \bigl(\partial _{t}\gamma +v_{\gamma ,\lambda }\bigr)\,d\sigma =0.  \label{eq:weak-curve}
\end{equation}%
Since $(\tau ,n)$ is an orthonormal frame, the definition~%
\eqref{eq:gauge-function} gives $\partial _{t}\gamma +v_{\gamma ,\lambda }=B\,n+g\,\tau $, where $B:=(\partial _{t}\gamma +v_{\gamma ,\lambda })\cdot n$. The chain rule
gives $\partial _{\sigma }[\phi (\gamma )]=|\partial _{\sigma }\gamma
|\,\nabla \phi (\gamma )\cdot \tau $. We integrate by parts on $\mathbb{T}$,
where periodicity removes the boundary terms, and obtain
\begin{equation*}
\int_{\mathbb{T}}\lambda \,g\,\nabla \phi (\gamma )\cdot \tau \,d\sigma
=\int_{\mathbb{T}}\frac{\lambda \,g}{|\partial _{\sigma }\gamma |}\,\partial
_{\sigma }[\phi (\gamma )]\,d\sigma =-\int_{\mathbb{T}}\partial _{\sigma }%
\Bigl(\frac{\lambda \,g}{|\partial _{\sigma }\gamma |}\Bigr)\phi (\gamma
)\,d\sigma .
\end{equation*}%
Write $A:=\partial _{t}\lambda -\partial _{\sigma }(\lambda g/|\partial _{\sigma }\gamma |)$ and $\partial _{n}\phi :=\nabla \phi \cdot n$. Then~%
\eqref{eq:weak-curve} reads
\begin{equation}
\int_{\mathbb{T}}A\,\phi (\gamma )\,d\sigma +\int_{\mathbb{T}}\lambda
\,B\,\partial _{n}\phi (\gamma )\,d\sigma =0.  \label{eq:weak-split}
\end{equation}%
If~\eqref{eq:normal-law} and~\eqref{eq:continuity-lambda} hold, then $B=0$ and $A=0$, so~\eqref{eq:weak-split} holds for every $\phi $.

Conversely, assume that~\eqref{eq:weak-split} holds for every $\phi $, and
fix $t$. Since $\gamma (\cdot ,t)$ is a smooth embedding of $\mathbb{T}$, the
tubular neighborhood theorem~\cite[Chapter~6]{Lee13} gives $\varepsilon >0$ such
that $\Psi (\sigma ,s):=\gamma (\sigma ,t)+s\,n(\sigma ,t)$ is a
diffeomorphism from $\mathbb{T}\times (-\varepsilon ,\varepsilon )$ onto an open
neighborhood $U$ of $\Gamma _{t}$. Let $a,b\in C^{\infty }(\mathbb{T})$, and
let $\chi \in C_{c}^{\infty }((-\varepsilon ,\varepsilon ))$ be equal to $1$ near $s=0$.
Define $\phi (\Psi (\sigma ,s)):=\chi (s)\bigl(a(\sigma )+s\,b(\sigma )\bigr)$
on $U$ and $\phi :=0$ outside $U$. The function $\phi $ is smooth on $U$ and
vanishes outside the compact set $\Psi (\mathbb{T}\times \mathrm{supp}\,\chi
)\subset U$, so $\phi \in C_{c}^{\infty }(\mathbb{R}^{2})$. At $s=0$ we have $%
\phi (\gamma (\sigma ,t))=a(\sigma )$. Since $\partial _{s}\Psi =n$, we also
have $\partial _{n}\phi (\gamma (\sigma ,t))=\partial _{s}[\phi (\Psi
(\sigma ,s))]_{s=0}=b(\sigma )$. Then~\eqref{eq:weak-split} becomes
\begin{equation*}
\int_{\mathbb{T}}A\,a\,d\sigma +\int_{\mathbb{T}}\lambda \,B\,b\,d\sigma =0\qquad \text{for all }a,b\in C^{\infty }(\mathbb{T}).
\end{equation*}%
The choice $b=0$ gives $\int_{\mathbb{T}}A\,a\,d\sigma =0$ for every smooth $%
a$. Since $A$ is continuous, $A=0$, which is~\eqref{eq:continuity-lambda}.
The choice $a=0$ gives $\lambda \,B=0$ in the same way. Since $%
\lambda >0$, we get $B=0$, which is~\eqref{eq:normal-law}.
\end{proof}\medskip

\noindent
As a consequence of the above result, we have the following complementary characterization of weak solutions in the Lagrangian formulation.
\begin{corollary}[Lagrangian labels]
\label{cor:lagrangian-labels}Assume the hypotheses of
Proposition~\ref{prop:euler-lagrange}.
\begin{enumerate}
\item[(i)] If $\partial _{t}\gamma =-v_{\gamma ,\lambda }$, then $\mu _{t}$
solves~\eqref{eq:transport} weakly if and only if $\partial _{t}\lambda =0$.
\item[(ii)] If $\mu _{t}$ solves~\eqref{eq:transport} weakly, then there exists a relabeling of $\Gamma _{t}$ that depends smoothly on $t$, is the identity at $t=0$, and gives a pair $(\tilde{\gamma},\tilde{\lambda})$ that describes the same measure $\mu _{t}$ and satisfies the Lagrangian laws~\eqref{eq:lagrangian-intro}.
\end{enumerate}
\end{corollary}
\begin{proof}
(i) If $\partial _{t}\gamma =-v_{\gamma ,\lambda }$, then $g=0$ by~\eqref{eq:gauge-function},
and~\eqref{eq:normal-law} holds. By Proposition~\ref{prop:euler-lagrange},
$\mu _{t}$ is a weak solution if and only if~\eqref{eq:continuity-lambda}
holds with $g=0$, that is, $\partial _{t}\lambda =0$.

(ii) We move each label along $\Gamma _{t}$ so as to cancel its tangential velocity $g$. Let $\psi $ solve $\partial _{t}\psi (\sigma ,t)=-g(\psi ,t)/|\partial
_{\sigma }\gamma (\psi ,t)|$ with $\psi (\sigma ,0)=\sigma $. The right-hand
side is a smooth time-dependent vector field on the compact manifold
$\mathbb{T}$. Its flow therefore exists on $[0,T)$, and each $\psi (\cdot ,t)$
is a diffeomorphism isotopic to the identity~\cite[Chapter~9]{Lee13}, hence
orientation preserving, with $\partial _{\sigma }\psi >0$. We set $\tilde{\gamma}(\sigma ,t):=\gamma (\psi (\sigma ,t),t)$ and $\tilde{\lambda}(\sigma ,t):=\lambda (\psi (\sigma ,t),t)\,\partial _{\sigma }\psi (\sigma ,t)$. By~\eqref{eq:normal-law} and~\eqref{eq:gauge-function},
$\partial _{t}\gamma =-v_{\gamma ,\lambda }+g\,\tau $. The chain rule gives
\begin{align*}
\partial _{t}\tilde{\gamma}& =(\partial _{t}\gamma )(\psi ,t)+(\partial
_{\sigma }\gamma )(\psi ,t)\,\partial _{t}\psi \\
& =-v_{\gamma ,\lambda }(\psi ,t)+g(\psi ,t)\,\tau (\psi ,t)-g(\psi
,t)\,\tau (\psi ,t)=-v_{\gamma ,\lambda }(\psi ,t).
\end{align*}%
The rule $\tilde{\lambda}\,d\sigma =\lambda \,d\psi $ is the transformation
rule of Section~\ref{subsec:contour-law}, so $(\tilde{\gamma},\tilde{\lambda})$
describes the same measure $\mu _{t}$. The velocity at a point of the curve
depends only on the measure and on the point, so $v_{\tilde{\gamma},\tilde{%
\lambda}}(\sigma ,t)=v_{\gamma ,\lambda }(\psi (\sigma ,t),t)$. Here we use
that the principal value does not change under a smooth change of parameter.
Indeed, the symmetric interval of length $2\varepsilon $ removed in one
parameter corresponds in the other parameter to an interval whose two halves
differ in length by $O(\varepsilon ^{2})$. On this difference the kernel is
$O(1/\varepsilon )$, so the two principal values differ by $O(\varepsilon )$.
Hence $\partial _{t}\tilde{\gamma}=-v_{\tilde{\gamma},\tilde{\lambda}}$. The
pair $(\tilde{\gamma},\tilde{\lambda})$ satisfies the hypotheses of
Proposition~\ref{prop:euler-lagrange}, since $\partial _{\sigma }\psi >0$, and
part~(i) gives $\partial _{t}\tilde{\lambda}=0$.
\end{proof}

\begin{remark}[The tangential gauge]
\label{rem:tangential-gauge}We call $g$ the tangential gauge. By~\eqref{eq:normal-law}, every curve-supported solution moves by $\partial _{t}\gamma =-v_{\gamma ,\lambda }+g\,\tau $. The normal velocity is fixed, and the term $g\,\tau $ only moves the labels along $\Gamma _{t}$. With the density per unit arc length $\rho =\lambda /|\partial _{\sigma }\gamma |$, the law~\eqref{eq:continuity-lambda} reads $\partial _{t}\lambda =\partial _{\sigma }(\rho \,g)$, so $g$ also decides how the line density changes. Relabeling the curve changes $g$ but not $\mu _{t}$. A relabeling that depends on $t$ slides each label along $\Gamma _{t}$, and the speed of this sliding is added to $g$. In the proof of Corollary~\ref{cor:lagrangian-labels}(ii) the labels slide with speed $-g$, which gives $g=0$, the Lagrangian labels. With the choice $g=v_{\gamma ,\lambda }\cdot \tau $, the curve moves only in the normal direction. All choices of $g$ describe the same measures $\mu _{t}$.
\end{remark}

\begin{remark}[Vortex sheets and patches]
These results have classical counterparts for the Euler equations of an ideal fluid. For a vortex sheet, the two one-sided velocities differ only in their tangential components; advecting the sheet by their average, which is given by a principal-value integral, makes the circulation a Lagrangian parameter, the analogue of the constant-density parameter of Section~\ref{subsec:gradient-flow}, and in this parameter the sheet obeys the Birkhoff--Rott equation~\cite[Section~9.2]{MB02}. For a vortex patch with smooth boundary, the weak solutions correspond exactly to the solutions of the contour dynamics equation for the boundary, written in a Lagrangian parametrization~\cite[Section~8.3.1, Proposition~8.6]{MB02}.
\end{remark}

We use Lagrangian labels from now on, because in them the contour law is a gradient flow.

\subsection{The contour law as a gradient flow}

\label{subsec:gradient-flow}

By~\eqref{eq:gradX-final} and~\eqref{eq:velocity}, the gradient of the energy $E_{\lambda }(\gamma )$ of Section~\ref{sec:energy-curves} with respect to $\gamma $, at fixed $\lambda $, is
\begin{equation*}
\nabla _{\gamma }E_{\lambda }(\gamma )(\sigma ,t)=\lambda (\sigma
)\,v_{\gamma ,\lambda }(\sigma ,t),
\end{equation*}%
so that the Lagrangian evolution \eqref{eq:lagrangian-intro} takes the form
\begin{equation}
\lambda (\sigma )\,\partial _{t}\gamma (\sigma ,t)=-\nabla _{\gamma
}E_{\lambda }(\gamma )(\sigma ,t).  \label{eq:gradientflow}
\end{equation}

A \emph{Lagrangian equilibrium} is a pair $(\gamma,\lambda)$ with $%
\partial_t\gamma = 0$ and $\partial_t\lambda = 0$. Inserting into %
\eqref{eq:gradientflow} yields the equilibrium condition
\begin{equation}  \label{eq:LagrangianEquilibrium}
\nabla_{\gamma}E_{\lambda}(\gamma)(\sigma) = 0 \qquad \forall\,\sigma \in
\mathbb{T},
\end{equation}
which coincides with the critical-point equation derived in %
\eqref{eq:gradX-final}. We call $(\gamma _{0},\lambda _{0})$ an equilibrium when the induced velocity vanishes on $\Gamma _{0}$, $v_{\gamma _{0},\lambda _{0}}\equiv 0$.

Equation~\eqref{eq:gradientflow} is a \emph{weighted} gradient flow of $%
E_{\lambda }$ with respect to $\gamma $, where the weight at each label $%
\sigma $ is $\lambda (\sigma )$. Since $\lambda >0$, it is the gradient flow of $E_{\lambda }$ for the metric $\langle h_{1},h_{2}\rangle _{\lambda }=\int_{\mathbb{T}}\lambda \,h_{1}\cdot h_{2}\,d\sigma $. In a
general parametrization the weight varies along the curve.

\paragraph{Constant-density parametrization.}

A parametrization in which $\lambda $ is constant removes the weight
in~\eqref{eq:gradientflow}, so the contour law becomes an unweighted gradient flow.

Such a parametrization always exists. Let $\Gamma $ have length $L$ and arc-length parameter $s\in \lbrack 0,L)$, let $\rho >0$ be its smooth density per unit arc length, and let $M=\int_{0}^{L}\rho \,ds$ be the total mass. Since $\rho >0$, the cumulative mass $\Phi (s)=\frac{2\pi }{M}\int_{0}^{s}\rho (s^{\prime })\,ds^{\prime }$ is smooth and strictly increasing from $0$ to $2\pi $, so $\sigma =\Phi (s)$ is a new parameter on $\mathbb{T}$ and $\tilde{\gamma}=\gamma \circ \Phi ^{-1}$ is a reparametrization of $\Gamma $. As the mass of every arc is preserved, $\tilde{\lambda}(\sigma )\,d\sigma =\rho (s)\,ds$, and therefore $\tilde{\lambda}=\rho /\Phi ^{\prime }=M/2\pi $ is constant. We call $\tilde{\gamma}$ the \emph{constant-density parametrization} of $\Gamma $: every parameter interval of length $\Delta \sigma $ carries the mass $(M/2\pi )\,\Delta \sigma $, wherever it lies on $\Gamma $.

We apply this reparametrization once, at $t=0$, and then use Lagrangian labels. By~\eqref{eq:lagrangian-intro} the line density does not change in time, and it equals $\tilde{\lambda}=M/2\pi $ for all $\sigma $ and all $t$; only the curve $\tilde{\gamma}(\sigma ,t)$ evolves.
In the constant-density parametrization, the constant $\tilde{\lambda}=M/2\pi $
factors out of each term of the energy, which becomes
\begin{equation}
E_{\tilde{\lambda}}(\tilde{\gamma})=\frac{1}{2}\left( \frac{M}{2\pi }\right)
^{2}\iint_{\mathbb{T}\times \mathbb{T}}G\bigl(\tilde{\gamma}(\sigma ,t),%
\tilde{\gamma}(\sigma ^{\prime },t)\bigr)\,d\sigma \,d\sigma ^{\prime
}+\omega \,\frac{M}{2\pi }\int_{\mathbb{T}}V\bigl(\tilde{\gamma}(\sigma ,t)%
\bigr)\,d\sigma. \label{eq:energy-uniform}
\end{equation}%
Note that the interaction term is proportional to $\tilde{\lambda}^{2}$ and the confinement term is proportional to $\tilde{\lambda}$. The weighted gradient flow %
\eqref{eq:gradientflow} reduces to
\begin{equation}
\partial _{t}\tilde{\gamma}(\sigma ,t)=-\frac{1}{\tilde{\lambda}}\,\nabla _{%
\tilde{\gamma}}E_{\tilde{\lambda}}(\tilde{\gamma})(\sigma ,t),
\label{eq:gradientflow-uniform}
\end{equation}%
which is a true gradient flow of $E_{\tilde{\lambda}}$ with respect to $%
\tilde{\gamma}$, with a uniform (label-independent) metric weight $\tilde{\lambda}$.

\section{Instability of equilibrium curves}

\label{sec:geometric}

In this section, we study the stability of an equilibrium, that is, of a critical point of the energy (Section~\ref{subsec:gradient-flow}). In Lagrangian labels the contour law is a gradient flow, so its linearization is given by the second variation of the energy expressed in these coordinates. Section~\ref{subsec:linearization-circle} computes the linearization exactly on the circle with quadratic confinement, and Section~\ref{subsec:energy-instability} proves, directly from the energy, that every smooth equilibrium curve is unstable, in the sense that it is not a local minimizer, irrespective of the choice of confining potential.

Let $\gamma _{0}$ be a smooth embedding of $\mathbb{T}$ and let $\lambda _{0}>0$ be smooth, with $v_{\gamma _{0},\lambda _{0}}\equiv 0$, so that $(\gamma _{0},\lambda _{0})$ is an equilibrium. We write $\tau _{0}$, $n_{0}$ and $\rho _{0}=\lambda _{0}/|\partial _{\sigma }\gamma _{0}|$ for the unit tangent, the outward unit normal and the density per unit arc length of $\Gamma _{0}$. By Corollary~\ref{cor:lagrangian-labels}, the line density stays equal to $\lambda _{0}$ in Lagrangian labels, so we write $E(\gamma )=E_{\lambda _{0}}(\gamma )$ and the contour law~\eqref{eq:gradientflow} reads $\lambda _{0}\,\partial _{t}\gamma =-\nabla _{\gamma }E(\gamma )$. We perturb the curve,
\begin{equation}
\gamma (\sigma ,t)=\gamma _{0}(\sigma )+\varepsilon \,h(\sigma ,t),\qquad h=\eta \,n_{0}+\xi \,\tau _{0},  \label{eq:linearization-ansatz}
\end{equation}%
with a normal displacement $\eta $ and a tangential displacement $\xi $.

The energy depends smoothly on the curve (Section~\ref{subsec:first-variation}), and its gradient $\nabla _{\gamma }E(\gamma )=\lambda _{0}v_{\gamma ,\lambda _{0}}$, given by~\eqref{eq:gradX-final}, vanishes at $\gamma _{0}$. So, with $\nabla ^{2}E$ denoting the Hessian of $E$ at $\gamma _{0}$, the energy has no term of order $\varepsilon $, and we may write
\begin{equation}
  \begin{gathered}
    E(\gamma _{0}+\varepsilon h)=E(\gamma _{0})+\varepsilon ^{2}E_{2}(h)+O(\varepsilon ^{3}),\\
    \text{where}\qquad E_{2}(h)=\frac{1}{2}\int_{\mathbb{T}}h\cdot \nabla ^{2}E\,h\,d\sigma.
    \end{gathered}\label{eq:second-variation}
\end{equation}%
Differentiating $\nabla _{\gamma }E(\gamma )=\lambda _{0}v_{\gamma ,\lambda _{0}}$ gives $\nabla ^{2}E\,h=\lambda _{0}\,v_{1}$, where $v_{1}$ is the term of order $\varepsilon $ of $v_{\gamma _{0}+\varepsilon h,\lambda _{0}}$. So the terms of order $\varepsilon $ of the contour law give, formally, the linear equation
\begin{equation}
\lambda _{0}\,\partial _{t}h=-\nabla ^{2}E\,h=-\lambda _{0}\,v_{1}.  \label{eq:lin-hessian}
\end{equation}%

\subsection{Exact growth rates on the circle}

\label{subsec:linearization-circle}

We now take the circle. We identify $\mathbb{R}^{2}$ with $\mathbb{C}$ and use the quadratic trap~\eqref{eq:V-quadratic}. For $\omega >0$ and $M>0$ let
\begin{equation}
r^{2}=\frac{M}{4\pi \omega },\qquad \gamma _{0}(\sigma )=r\,e^{i\sigma },\qquad \lambda _{0}=\frac{M}{2\pi }.  \label{eq:circle-data}
\end{equation}%
This is the circle of Section~\ref{sec:circle-example}. The angle $\sigma $ is the parameter, and $|\partial _{\sigma }\gamma _{0}|=r$. The relation~\eqref{eq:circle-data} reads $\omega =\lambda _{0}/(2r^{2})$. The unit normal and the unit tangent are $n_{0}=e^{i\sigma }$ and $\tau _{0}=ie^{i\sigma }$, and $\eta $ and $\xi $ are real. We write $f(\sigma )=\sum_{k\in \mathbb{Z}}\widehat{f}_{k}e^{ik\sigma }$ for the corresponding Fourier series.
\begin{theorem}[Linear spectrum of the circle]
\label{thm:circle-spectrum}The pair~\eqref{eq:circle-data} is an equilibrium, $v_{\gamma _{0},\lambda _{0}}\equiv 0$. For every smooth displacement $h=\eta \,n_{0}+\xi \,\tau _{0}$,
\begin{equation}
E_{2}(h)=\pi \omega \lambda _{0}\sum_{k\in \mathbb{Z}}\Bigl(-(|k|-2)\,|\widehat{\eta }_{k}|^{2}+|k|\,|\widehat{\xi }_{k}|^{2}\Bigr).  \label{eq:circle-Q}
\end{equation}%
The linearized flow~\eqref{eq:lin-hessian} decouples into
\begin{equation}
\partial _{t}\eta =\omega \bigl(|D|-2\bigr)\eta ,\qquad \partial _{t}\xi =-\omega \,|D|\xi .  \label{eq:circle-linear}
\end{equation}%
Hence, for every $k\in \mathbb{Z}$,
\begin{equation}
\widehat{\eta }_{k}(t)=e^{\omega (|k|-2)t}\,\widehat{\eta }_{k}(0),\qquad \widehat{\xi }_{k}(t)=e^{-\omega |k|t}\,\widehat{\xi }_{k}(0).  \label{eq:circle-modes}
\end{equation}%
The normal displacement grows in the modes $|k|\geq 3$, is neutral in the modes $|k|=2$, and decays in the modes $|k|\leq 1$. The tangential displacement decays in every mode $k\neq 0$.
\end{theorem}
\begin{proof}
\emph{The velocity in complex form.} We identify $\mathbb{R}^{2}$ with $\mathbb{C}$, so that $x/|x|^{2}=1/\overline{x}$. Since $\nabla _{x}G(x,y)=-(x-y)/(2\pi |x-y|^{2})$ and $\nabla V(x)=x$, the velocity~\eqref{eq:velocity} of a curve $\gamma $ with the density $\lambda _{0}$ is
\begin{equation*}
v_{\gamma ,\lambda _{0}}(\sigma )=-\frac{\lambda _{0}}{2\pi }\,\mathrm{p.v.}\int_{\mathbb{T}}\frac{d\sigma ^{\prime }}{\overline{\delta \gamma }}+\omega \gamma (\sigma ),
\end{equation*}%
where $\delta \gamma =\gamma (\sigma )-\gamma (\sigma ^{\prime })$ is the chord, in the notation of Definition~\ref{def:delta}.\smallskip

\noindent
\emph{The kernel on the circle.} Since $\gamma _{0}=r\,e^{i\sigma }$,
\begin{equation*}
\overline{\delta \gamma _{0}}=r\,e^{-i\sigma }(1-e^{i\delta \sigma })=-2ir\,e^{-i\sigma }e^{i\delta \sigma /2}\sin \frac{\delta \sigma }{2}.
\end{equation*}%
Since $e^{-i\delta \sigma /2}/\sin \frac{\delta \sigma }{2}=\cot \frac{\delta \sigma }{2}-i$ and $\tau _{0}=i\,n_{0}$, this gives
\begin{equation*}
\frac{1}{\overline{\delta \gamma _{0}}}=\frac{1}{2r}\Bigl(n_{0}(\sigma )+\cot \frac{\delta \sigma }{2}\,\tau _{0}(\sigma )\Bigr).
\end{equation*}%
Since the cotangent is odd, $\mathrm{p.v.}\int_{\mathbb{T}}\cot \frac{\delta \sigma }{2}\,d\sigma ^{\prime }=0$.\smallskip

\noindent
\emph{The equilibrium.} By the formula for $1/\overline{\delta \gamma _{0}}$ in the previous step, the cotangent term has principal value zero and the term $n_{0}(\sigma )/(2r)$ does not depend on $\sigma ^{\prime }$. So the principal value of $\int_{\mathbb{T}}d\sigma ^{\prime }/\overline{\delta \gamma _{0}}$ is $2\pi \cdot n_{0}(\sigma )/(2r)=\frac{\pi }{r}\,n_{0}(\sigma )$. Hence
\begin{equation*}
v_{\gamma _{0},\lambda _{0}}=-\frac{\lambda _{0}}{2r}\,n_{0}+\omega r\,n_{0}=0,
\end{equation*}%
because $\omega =\lambda _{0}/(2r^{2})$. So the pair is an equilibrium, which confirms~\eqref{eq:rstar} without restricting to centered circles.\smallskip

\noindent
\emph{The linearization.} The chord of the perturbed curve $\gamma _{0}+\varepsilon h$ is $\delta \gamma _{0}+\varepsilon \,\delta h$, and $1/\overline{\delta \gamma _{0}+\varepsilon \,\delta h}=1/\overline{\delta \gamma _{0}}-\varepsilon \,\overline{\delta h}/\overline{\delta \gamma _{0}}^{2}+O(\varepsilon ^{2})$. So the term of order $\varepsilon $ of $v_{\gamma _{0}+\varepsilon h,\lambda _{0}}$ is
\begin{equation*}
v_{1}=\frac{\lambda _{0}}{2\pi }\,\overline{\mathcal{I}}+\omega h,\qquad \mathcal{I}=\mathrm{p.v.}\int_{\mathbb{T}}\frac{\delta h}{(\delta \gamma _{0})^{2}}\,d\sigma ^{\prime }.
\end{equation*}

\noindent
\emph{The integral.} Write $h=u\,e^{i\sigma }$ with $u=\eta +i\xi $. Then, by the formula for $\overline{\delta \gamma _{0}}$ above,
\begin{equation*}
\delta h=e^{i\sigma }\bigl(u(\sigma )-u(\sigma ^{\prime })e^{-i\delta \sigma }\bigr),\qquad (\delta \gamma _{0})^{2}=-4r^{2}e^{2i\sigma }e^{-i\delta \sigma }\sin ^{2}\frac{\delta \sigma }{2}.
\end{equation*}%
Hence $\delta h/(\delta \gamma _{0})^{2}=-\frac{1}{4r^{2}}\,e^{-i\sigma }B$, where
\begin{equation*}
B=\frac{u(\sigma )\,e^{i\delta \sigma }-u(\sigma ^{\prime })}{\sin ^{2}\frac{\delta \sigma }{2}}=\frac{\delta u+u(\sigma )\,(e^{i\delta \sigma }-1)}{\sin ^{2}\frac{\delta \sigma }{2}}=\frac{\delta u}{\sin ^{2}\frac{\delta \sigma }{2}}+2iu(\sigma )\cot \frac{\delta \sigma }{2}-2u(\sigma ),
\end{equation*}%
since $e^{i\delta \sigma }-1=2i\sin \frac{\delta \sigma }{2}\cos \frac{\delta \sigma }{2}-2\sin ^{2}\frac{\delta \sigma }{2}$. By Lemma~\ref{lem:piD}, the principal value of $\int_{\mathbb{T}}\delta u/\sin ^{2}\frac{\delta \sigma }{2}\,d\sigma ^{\prime }$ is $4\pi \,|D|u$. The term $2iu(\sigma )\cot \frac{\delta \sigma }{2}$ has principal value zero, and $\int_{\mathbb{T}}2u(\sigma )\,d\sigma ^{\prime }=4\pi u(\sigma )$. Hence
\begin{equation*}
\mathcal{I}=-\frac{e^{-i\sigma }}{4r^{2}}\bigl(4\pi \,|D|u-4\pi u\bigr)=-\frac{\pi }{r^{2}}\,e^{-i\sigma }\,(|D|-1)u.
\end{equation*}%

\noindent
\emph{The Hessian.} Since $\lambda _{0}/(2r^{2})=\omega $ and $|D|$ maps real functions to real functions,
\begin{equation*}
v_{1}=-\omega \,e^{i\sigma }\,\overline{(|D|-1)u}+\omega \,e^{i\sigma }u=\omega \,e^{i\sigma }\bigl((2-|D|)\eta +i\,|D|\xi \bigr).
\end{equation*}%
Since $n_{0}=e^{i\sigma }$ and $\tau _{0}=ie^{i\sigma }$, this gives
\begin{equation*}
\nabla ^{2}E\,h=\lambda _{0}v_{1}=-\omega \lambda _{0}\bigl((|D|-2)\eta \bigr)\,n_{0}+\omega \lambda _{0}\bigl(|D|\xi \bigr)\,\tau _{0}.
\end{equation*}%

\noindent
\emph{Fourier modes.} By~\eqref{eq:second-variation},
\begin{equation*}
E_{2}(h)=\frac{1}{2}\int_{\mathbb{T}}h\cdot \nabla ^{2}E\,h\,d\sigma =\frac{\omega \lambda _{0}}{2}\int_{\mathbb{T}}\Bigl(-\eta \,(|D|-2)\eta +\xi \,|D|\xi \Bigr)\,d\sigma .
\end{equation*}%
By~\eqref{eq:piD-identity}, $|D|$ multiplies $e^{ik\sigma }$ by $|k|$. For real $f$ and $g$, Parseval's identity reads $\int_{\mathbb{T}}f\,g\,d\sigma =2\pi \sum_{k}\widehat{f}_{k}\,\overline{\widehat{g}_{k}}$, and it turns the last display into~\eqref{eq:circle-Q}. By~\eqref{eq:lin-hessian}, $\partial _{t}h=-v_{1}$, and its components along $n_{0}$ and $\tau _{0}$, which do not depend on $t$, give~\eqref{eq:circle-linear}. In Fourier modes this reads $\partial _{t}\widehat{\eta }_{k}=\omega (|k|-2)\,\widehat{\eta }_{k}$ and $\partial _{t}\widehat{\xi }_{k}=-\omega |k|\,\widehat{\xi }_{k}$, which gives~\eqref{eq:circle-modes}.
\end{proof}\medskip

The modes with $|k|\leq 1$ have a direct meaning. The normal mode $k=0$ is a dilation, which decays like $e^{-2\omega t}$. A translation has modes $|k|=1$ in both $\eta $ and $\xi $, and both decay like $e^{-\omega t}$, because the trap pulls the circle back to the origin. The tangential mode $k=0$ is a rigid rotation, which only relabels the circle and is neutral. By~\eqref{eq:circle-Q}, $E_{2}$ is negative on normal displacements made of modes with $|k|\geq 3$, so the circle is not a local minimizer of $E$. Since the growth rates $\omega (|k|-2)$ are unbounded in $k$, the linearized evolution on the circle is ill-posed in the sense of Hadamard, as for the Kelvin--Helmholtz instability of vortex sheets~\cite{Moore79,Krasny86}, where the $k$th Fourier mode of the linearization about a flat sheet has a component that grows like $e^{|k|t/2}$~\cite[Section~9.3]{MB02}.

\paragraph{Comparison with the regular polygon.}
The spectrum of the circle is the limit of that of the regular $N$-gon of Section~\ref{sec:circle-example}. The linearization of the polygon is classical~\cite{CS00,GI12}. In the discrete Fourier mode $k\in (-N/2,N/2]\cap \mathbb{Z}$ of the vertices, the Hessian of $E_{N}/N$ at the $N$-gon is diagonal in the radial and tangential displacements, and its two entries are $\omega (2-d_{k})$ for the normal and $\omega d_{k}$ for the tangential displacement, where
\begin{equation*}
d_{k}=\frac{|k|(N-|k|)}{N-1}=|k|-\frac{|k|(|k|-1)}{N-1}
\end{equation*}
takes the place of $|k|$ in Theorem~\ref{thm:circle-spectrum}. For each fixed $k$, $d_{k}$ tends to $|k|$ as $N\rightarrow \infty $, so the spectrum of the polygon converges, mode by mode, to that of the circle.

This convergence is not uniform in $k$, and for small $N$ it changes the stability. The vortex spacing is a length scale that the continuum does not have. For the gradient dynamics, the regular $N$-gon is a nondegenerate local minimizer of the energy, up to rotation, for $N\leq 6$, it is degenerate for $N=7$, and it is unstable for $N\geq 8$. Small polygons are therefore stable, although the circle that they approximate is not. The rotating $N$-gon is nonlinearly stable for $N\leq 7$, the case $N=7$ being proved in~\cite{KY02}, and unstable for $N\geq 8$.

\subsection{Instability of the energy}

\label{subsec:energy-instability}

We now describe, formally, the linear equation~\eqref{eq:lin-hessian} on a general equilibrium curve. In Lagrangian labels the unknowns are the two components of the displacement $h=\eta \,n_{0}+\xi \,\tau _{0}$ in~\eqref{eq:linearization-ansatz}. The normal displacement $\eta $ moves the curve and changes its shape. The tangential displacement $\xi $ moves the labels along $\Gamma _{0}$. To first order, it leaves the curve unchanged and only redistributes its mass along it, so that in the original parametrization the line density changes by $-\varepsilon \,\partial _{\sigma }(\rho _{0}\xi )$. The decomposition~\eqref{eq:grad-final-general} for the normal part, and the fact that the velocity induced on $\Gamma _{0}$ by a density is, to leading order, tangential, lead formally, after some computations, to the system
\begin{equation}
\begin{aligned}
\partial _{t}\eta & =\frac{\rho _{0}}{2|\partial _{\sigma }\gamma _{0}|}\,|D|\eta +\mathcal{B}_{1}(\eta ,\xi ), \\
\partial _{t}\xi & =-\frac{\rho _{0}}{2|\partial _{\sigma }\gamma _{0}|}\,|D|\xi +\mathcal{B}_{2}(\eta ,\xi ).
\end{aligned}
\label{eq:formal-system}
\end{equation}%
Here, $\mathcal{B}_{1}$ and $\mathcal{B}_{2}$ collect the remaining lower-order terms. On the circle the coefficient in front of $|D|$ is $\omega $, $\mathcal{B}_{1}(\eta ,\xi )=-2\omega \eta $ and $\mathcal{B}_{2}=0$, which is exactly~\eqref{eq:circle-linear}.

On a general curve, the structure of the kernel of the remainder $\mathcal{R}_{\lambda _{0}}$ in~\eqref{eq:R-def} indicates that $\mathcal{B}_{1}$ and $\mathcal{B}_{2}$ are of lower order. Normal perturbations of high wavenumber $k$ would then grow at rates of order $|k|$, while the leading part of the tangential equation is dissipative. Making this rigorous requires estimating the contribution of $\mathcal{R}_{\lambda _{0}}$ to the second variation for every displacement, which we do not do. Instead, to prove instability, it is enough to estimate it in the direction that is unstable on the circle and in~\eqref{eq:formal-system}, the normal displacement $\cos (k\sigma )\,n_{0}$ with $k$ large. This uses only the quadratic form $E_{2}$ of~\eqref{eq:second-variation}, and it shows that no smooth equilibrium curve is a local minimizer.

\begin{theorem}[Second variation on a general curve]
\label{thm:energy-instability}Let $h_{k}=\cos (k\sigma )\,n_{0}$ for $k\in \mathbb{N}$, where $(\gamma _{0},\lambda _{0})$ is an equilibrium as at the beginning of this section. Then, as $k\rightarrow \infty $,
\begin{equation*}
E_{2}(h_{k})=-\frac{k}{8}\int_{\mathbb{T}}\rho _{0}^{2}\,d\sigma +O(1).
\end{equation*}%
In particular, $E_{2}(h_{k})<0$ for $k$ large, so that $E(\gamma _{0}+\varepsilon h_{k})<E(\gamma _{0})$ for every sufficiently small $\varepsilon \neq 0$, and hence $\gamma _{0}$ is not a local minimizer of $E=E_{\lambda _{0}}$.
\end{theorem}

\begin{proof}
Let $h=h_{k}=\eta \,n_{0}$ with $\eta =\cos k\sigma $. In Lagrangian labels the density does not move, so $\lambda =\lambda _{0}$ throughout the proof, and $E=E_{\lambda _{0}}$ depends only on the curve. By~\eqref{eq:second-variation}, $E_{2}(h)=\frac{1}{2}\int_{\mathbb{T}}h\cdot \nabla ^{2}E\,h\,d\sigma $, where $\nabla ^{2}E\,h$ is the derivative of the gradient~\eqref{eq:grad-final-general} at $\gamma _{0}$ in the direction $h$. We differentiate the three terms of~\eqref{eq:grad-final-general}, namely the singular part, the trap and the remainder.
\smallskip

\noindent\emph{Orders.} Below, $O(1)$ denotes a quantity bounded by a constant that does not depend on $k$. For a function of $(\sigma ,\sigma ^{\prime })$, which may depend on $k$, we write $\mathcal{O}(|\delta \sigma |^{m})$ for a bound by $C|\delta \sigma |^{m}$, where $\delta \sigma $ is taken in $(-\pi ,\pi ]$ and the constant $C$ does not depend on $\sigma $, $\sigma ^{\prime }$ or $k$. For instance, $\delta f=\mathcal{O}(|\delta \sigma |)$ for every smooth function $f$ that does not depend on $k$, $1/(4\sin ^{2}\frac{\delta \sigma }{2})=\mathcal{O}(|\delta \sigma |^{-2})$, and $\eta (\sigma )$, $\eta (\sigma ^{\prime })$ and $\delta \eta $ are $\mathcal{O}(1)$.\smallskip

\noindent
\emph{The singular part.} The singular part of~\eqref{eq:grad-final-general} is $-\frac{\lambda _{0}^{2}}{2|\partial _{\sigma }\gamma |^{2}}\,|D|\gamma $. Since $\partial _{\sigma }\gamma _{0}\cdot n_{0}=0$, the derivative of $|\partial _{\sigma }\gamma |^{2}$ in the direction $h$ is $2\eta \,\partial _{\sigma }\gamma _{0}\cdot \partial _{\sigma }n_{0}$, which contains no derivative of $\eta $. So the derivative of the singular part is $-\frac{1}{2}\rho _{0}^{2}\,|D|h$ plus $\eta $ times a smooth vector field that does not depend on $k$, and $\eta $ times this field contributes $O(1)$ to $E_{2}(h)$. Since $h=\eta \,n_{0}$, the contribution of $-\frac{1}{2}\rho _{0}^{2}\,|D|h$ is $-\frac{1}{4}\int_{\mathbb{T}}\rho _{0}^{2}\,\eta \,n_{0}\cdot |D|h\,d\sigma $.

By Lemma~\ref{lem:piD}, $\pi \,n_{0}\cdot |D|h$ is the principal value of $\int_{\mathbb{T}}n_{0}(\sigma )\cdot \delta h/(4\sin ^{2}\frac{\delta \sigma }{2})\,d\sigma ^{\prime }$. Since $n_{0}$ is a smooth unit vector field,
\begin{equation}
n_{0}(\sigma )\cdot \delta h=\delta \eta +\frac{1}{2}\eta (\sigma ^{\prime })\,|\delta n_{0}|^{2}.  \label{eq:thm8-ndh}
\end{equation}
In~\eqref{eq:thm8-ndh}, the term $\delta \eta $ gives $\pi |D|\eta =\pi k\eta $ by Lemma~\ref{lem:piD}, and the term $\frac{1}{2}\eta (\sigma ^{\prime })\,|\delta n_{0}|^{2}$ gives a bounded integrand, because it is $\mathcal{O}(|\delta \sigma |^{2})$. Therefore $n_{0}\cdot |D|h=k\eta +O(1)$, and the singular part contributes
\begin{equation*}
-\frac{k}{4}\int_{\mathbb{T}}\rho _{0}^{2}\cos ^{2}k\sigma \,d\sigma +O(1)=-\frac{k}{8}\int_{\mathbb{T}}\rho _{0}^{2}\,d\sigma +O(1),
\end{equation*}
because $\cos ^{2}k\sigma =\frac{1}{2}(1+\cos 2k\sigma )$ and $\int_{\mathbb{T}}\rho _{0}^{2}\cos 2k\sigma \,d\sigma =O(1/k)$ by one integration by parts.\smallskip

\noindent
\emph{The trap.} The derivative of the trap term is $\omega \lambda _{0}D^{2}V(\gamma _{0})h$, which contributes $O(1)$. Adding the contributions of the singular part and of the trap gives
\begin{equation}
E_{2}(h_{k})=-\frac{k}{8}\int_{\mathbb{T}}\rho _{0}^{2}\,d\sigma -\frac{1}{4\pi }\int_{\mathbb{T}}\lambda _{0}\,h\cdot D_{h}\mathcal{R}\,d\sigma +O(1),  \label{eq:thm8-E2}
\end{equation}
where $D_{h}\mathcal{R}$ is the derivative of $\mathcal{R}_{\lambda _{0}}$ at $\gamma _{0}$ in the direction $h$. It remains to show that the last integral in~\eqref{eq:thm8-E2} is bounded.\smallskip

\noindent
\emph{The remainder.} By~\eqref{eq:S-def} and~\eqref{eq:R-def}, $\mathcal{R}_{\lambda _{0}}(\gamma )(\sigma )=\int_{\mathbb{T}}\delta \gamma \,R\,d\sigma ^{\prime }$, with
\begin{equation*}
R=\frac{\lambda _{0}(\sigma ^{\prime })}{|\delta \gamma |^{2}}-\frac{\lambda _{0}(\sigma )}{|\partial _{\sigma }\gamma (\sigma )|^{2}\,4\sin ^{2}\frac{\delta \sigma }{2}},
\end{equation*}
and the integrand $\delta \gamma \,R$ is continuous by Proposition~\ref{prop:R-regular}. Let $R_{0}$ be $R$ at $\gamma _{0}$, and $D_{h}R$ the derivative of $R$ at $\gamma _{0}$ in the direction $h$. Since the chord of $\gamma _{0}+\varepsilon h$ is $\delta \gamma _{0}+\varepsilon \,\delta h$, the product rule gives
\begin{equation}
n_{0}(\sigma )\cdot \bigl(D_{h}\mathcal{R}\bigr)(\sigma )=\int_{\mathbb{T}}\bigl((n_{0}(\sigma )\cdot \delta h)\,R_{0}+(n_{0}(\sigma )\cdot \delta \gamma _{0})\,D_{h}R\bigr)\,d\sigma ^{\prime },  \label{eq:thm8-normal}
\end{equation}
where
\begin{equation}
D_{h}R=-\frac{2\lambda _{0}(\sigma ^{\prime })\,\delta \gamma _{0}\cdot \delta h}{|\delta \gamma _{0}|^{4}}+\frac{2\lambda _{0}(\sigma )\,\eta (\sigma )\,\partial _{\sigma }\gamma _{0}(\sigma )\cdot \partial _{\sigma }n_{0}(\sigma )}{|\partial _{\sigma }\gamma _{0}(\sigma )|^{4}\,4\sin ^{2}\frac{\delta \sigma }{2}}.  \label{eq:thm8-DhR}
\end{equation}
Since $h=\eta \,n_{0}$, the last integral in~\eqref{eq:thm8-E2} equals $\int_{\mathbb{T}}\lambda _{0}(\sigma )\eta (\sigma )\,n_{0}(\sigma )\cdot (D_{h}\mathcal{R})(\sigma )\,d\sigma $, and we bound the two terms in the integrand of~\eqref{eq:thm8-normal}.\smallskip

\noindent
\emph{The term with $D_{h}R$.} By Taylor's formula at $\sigma $ and at $\sigma ^{\prime }$, $\delta \gamma _{0}=\partial _{\sigma }\gamma _{0}(\sigma )\,\delta \sigma +\mathcal{O}(|\delta \sigma |^{2})=\partial _{\sigma }\gamma _{0}(\sigma ^{\prime })\,\delta \sigma +\mathcal{O}(|\delta \sigma |^{2})$. Because $\partial _{\sigma }\gamma _{0}$ is orthogonal to $n_{0}$ at the same point, $n_{0}(\sigma )\cdot \delta \gamma _{0}$ and $n_{0}(\sigma ^{\prime })\cdot \delta \gamma _{0}$ are $\mathcal{O}(|\delta \sigma |^{2})$, and so is $\delta \gamma _{0}\cdot \delta h=\eta (\sigma )\,n_{0}(\sigma )\cdot \delta \gamma _{0}-\eta (\sigma ^{\prime })\,n_{0}(\sigma ^{\prime })\cdot \delta \gamma _{0}$. Since $\gamma _{0}$ is a smooth embedding, $1/|\delta \gamma _{0}|=\mathcal{O}(|\delta \sigma |^{-1})$. Inserting these orders in~\eqref{eq:thm8-DhR} gives $D_{h}R=\mathcal{O}(|\delta \sigma |^{-2})$. Hence the term $(n_{0}(\sigma )\cdot \delta \gamma _{0})\,D_{h}R$ in~\eqref{eq:thm8-normal} is $\mathcal{O}(1)$.
\smallskip

\noindent
\emph{The term with $R_{0}$.} Since $\lambda _{0}(\sigma )^{2}/|\partial _{\sigma }\gamma _{0}(\sigma )|^{2}=\rho _{0}(\sigma )^{2}$, the definition of $R$ gives $\lambda _{0}(\sigma )R_{0}=\lambda _{0}(\sigma )\lambda _{0}(\sigma ^{\prime })/|\delta \gamma _{0}|^{2}-\rho _{0}(\sigma )^{2}/(4\sin ^{2}\frac{\delta \sigma }{2})$. We split this kernel into its symmetric and antisymmetric parts under the exchange of $\sigma $ and $\sigma ^{\prime }$. Only the numerator $\rho _{0}(\sigma )^{2}$ of the second term is not symmetric, and it is the sum of its symmetric part $\frac{1}{2}(\rho _{0}(\sigma )^{2}+\rho _{0}(\sigma ^{\prime })^{2})$ and its antisymmetric part $\frac{1}{2}\delta (\rho _{0}^{2})$. So the symmetric and antisymmetric parts of the kernel $\lambda _{0}(\sigma )R_{0}$ are the two terms of
\begin{equation}
\lambda _{0}(\sigma )R_{0}=\frac{F(\sigma ,\sigma ^{\prime })}{4\sin ^{2}\frac{\delta \sigma }{2}}-\frac{\delta (\rho _{0}^{2})}{8\sin ^{2}\frac{\delta \sigma }{2}},  \label{eq:thm8-R0}
\end{equation}
where
\begin{equation*}
F(\sigma ,\sigma ^{\prime })=\lambda _{0}(\sigma )\lambda _{0}(\sigma ^{\prime })\,\frac{4\sin ^{2}\frac{\delta \sigma }{2}}{|\delta \gamma _{0}|^{2}}-\frac{\rho _{0}(\sigma )^{2}+\rho _{0}(\sigma ^{\prime })^{2}}{2}.
\end{equation*}
The function $F$ is smooth and $F(\sigma ,\sigma )=0$. Let $\partial _{1}$ and $\partial _{2}$ denote the derivatives in the first and the second variable. Differentiating $F(\sigma ,\sigma )=0$ gives $\partial _{1}F(\sigma ,\sigma )+\partial _{2}F(\sigma ,\sigma )=0$, and differentiating the symmetry $F(\sigma ,\sigma ^{\prime })=F(\sigma ^{\prime },\sigma )$ in $\sigma $ and setting $\sigma ^{\prime }=\sigma $ gives $\partial _{1}F(\sigma ,\sigma )=\partial _{2}F(\sigma ,\sigma )$. So both first derivatives of $F$ vanish on the diagonal, and Taylor's formula in $\sigma ^{\prime }$ at $\sigma ^{\prime }=\sigma $ gives $F=\mathcal{O}(|\delta \sigma |^{2})$. Hence the symmetric part $F/(4\sin ^{2}\frac{\delta \sigma }{2})$ of the kernel is $\mathcal{O}(1)$.

Since $\delta (\rho _{0}^{2})=\mathcal{O}(|\delta \sigma |)$, the antisymmetric part of the kernel is $\mathcal{O}(|\delta \sigma |^{-1})$, and therefore $R_{0}=\mathcal{O}(|\delta \sigma |^{-1})$. The first term in the integrand of~\eqref{eq:thm8-normal} enters the last integral in~\eqref{eq:thm8-E2} multiplied by $\lambda _{0}(\sigma )\eta (\sigma )$. By~\eqref{eq:thm8-ndh},
\begin{equation}
\lambda _{0}(\sigma )\eta (\sigma )\bigl(n_{0}(\sigma )\cdot \delta h\bigr)R_{0}=\eta (\sigma )\,\delta \eta \,\lambda _{0}(\sigma )R_{0}+\tfrac{1}{2}\eta (\sigma )\eta (\sigma ^{\prime })\,|\delta n_{0}|^{2}\lambda _{0}(\sigma )R_{0}.  \label{eq:thm8-split}
\end{equation}
The last term of~\eqref{eq:thm8-split} is $\mathcal{O}(|\delta \sigma |^{2})\,\mathcal{O}(|\delta \sigma |^{-1})=\mathcal{O}(|\delta \sigma |)$. In the first term of~\eqref{eq:thm8-split} we insert~\eqref{eq:thm8-R0}. The symmetric part of the kernel gives $\eta (\sigma )\,\delta \eta \;\mathcal{O}(1)=\mathcal{O}(1)$, and the antisymmetric part of the kernel gives $-\frac{1}{2}$ times the integrand of the double integral in~\eqref{eq:thm8-main}. Together with the bound for the second term of~\eqref{eq:thm8-normal}, this gives
\begin{equation}
\int_{\mathbb{T}}\lambda _{0}\,h\cdot D_{h}\mathcal{R}\,d\sigma =-\frac{1}{2}\iint_{\mathbb{T}\times \mathbb{T}}\eta (\sigma )\,\delta \eta \,\frac{\delta (\rho _{0}^{2})}{4\sin ^{2}\frac{\delta \sigma }{2}}\,d\sigma ^{\prime }\,d\sigma +O(1).  \label{eq:thm8-main}
\end{equation}

\noindent
\emph{The double integral in~\eqref{eq:thm8-main}.} Exchanging the variables $\sigma $ and $\sigma ^{\prime }$ changes the signs of $\delta \eta $ and $\delta (\rho _{0}^{2})$. Since $\eta (\sigma )\,\delta \eta =\frac{1}{2}\delta (\eta ^{2})+\frac{1}{2}(\delta \eta )^{2}$, where $\delta (\eta ^{2})$ is antisymmetric and $(\delta \eta )^{2}$ is symmetric, the integrand of the double integral in~\eqref{eq:thm8-main} is the sum of its symmetric part $\frac{1}{2}\delta (\eta ^{2})\,\delta (\rho _{0}^{2})/(4\sin ^{2}\frac{\delta \sigma }{2})$ and its antisymmetric part $\frac{1}{2}(\delta \eta )^{2}\,\delta (\rho _{0}^{2})/(4\sin ^{2}\frac{\delta \sigma }{2})$. The antisymmetric part of the integrand is bounded for each $k$, and its double integral is zero, because exchanging the variables changes its sign but not the value of the integral. So the first term on the right of~\eqref{eq:thm8-main} equals
\begin{equation*}
-\frac{1}{4}\iint_{\mathbb{T}\times \mathbb{T}}\frac{\delta (\eta ^{2})\,\delta (\rho _{0}^{2})}{4\sin ^{2}\frac{\delta \sigma }{2}}\,d\sigma ^{\prime }\,d\sigma =-\frac{\pi }{2}\int_{\mathbb{T}}\rho _{0}^{2}\,|D|(\eta ^{2})\,d\sigma =-\frac{\pi k}{2}\int_{\mathbb{T}}\rho _{0}^{2}\cos 2k\sigma \,d\sigma .
\end{equation*}
The first equality is the identity $\iint \delta f\,\delta g/(4\sin ^{2}\frac{\delta \sigma }{2})\,d\sigma ^{\prime }\,d\sigma =2\pi \int f\,|D|g\,d\sigma $ for smooth $f$ and $g$, used with $f=\rho _{0}^{2}$ and $g=\eta ^{2}$. Both sides are bilinear, so by Fourier expansion it suffices to take $f=e^{im\sigma }$ and $g=e^{in\sigma }$. Then $\delta f\,\delta g=e^{i(m+n)\sigma }(1-e^{-im\,\delta \sigma })(1-e^{-in\,\delta \sigma })$, so both sides vanish unless $n=-m$. If $n=-m$, then $\delta f\,\delta g=4\sin ^{2}\frac{m\,\delta \sigma }{2}$, and the integral~\eqref{eq:fejer} in the proof of Lemma~\ref{lem:piD} shows that both sides equal $4\pi ^{2}|m|$. The second equality holds because $\eta ^{2}=\frac{1}{2}+\frac{1}{2}\cos 2k\sigma $, so $|D|(\eta ^{2})=k\cos 2k\sigma $. The integral $\int_{\mathbb{T}}\rho _{0}^{2}\cos 2k\sigma \,d\sigma $ is $O(1/k^{2})$ by two integrations by parts, so the right-hand side of the last display is $O(1/k)$.\smallskip

\noindent
\emph{Conclusion.} By~\eqref{eq:thm8-main}, the last integral in~\eqref{eq:thm8-E2} is therefore bounded, which proves the formula for $E_{2}(h_{k})$. Hence $E_{2}(h_{k})<0$ for $k$ large, and, by~\eqref{eq:second-variation}, $E(\gamma _{0}+\varepsilon h_{k})<E(\gamma _{0})$ for every small $\varepsilon \neq 0$.
\end{proof}\smallskip

On the circle, where $\rho _{0}=\lambda _{0}/r$ and $\omega =\lambda _{0}/(2r^{2})$, \eqref{eq:circle-Q} gives exactly $E_{2}(h_{k})=-\frac{\pi }{2}\,\omega \lambda _{0}(k-2)$, whose leading term is that of Theorem~\ref{thm:energy-instability}.

\section{Critical points on concentric circles}

\label{sec:numerics}

The single circle of Section~\ref{sec:circle-example} is the simplest
curve-supported critical point of the limit energy~\eqref{eq:E-intro}.
The energy also admits critical points concentrated on \emph{several}
concentric circles. We compute them explicitly for two circles in the
continuum, and then comment on configurations with more circles and on their discrete counterparts. Throughout we take unit mass and the trap strength $\omega _{\ast }=\frac{1}{2\pi }$, for which repulsion and confinement have the same coefficient. By rescaling, this is no loss of generality. With this choice, the energy~\eqref{eq:E-intro}
with the quadratic trap~\eqref{eq:V-quadratic} reads
\begin{equation}
E[\mu ]=-\frac{1}{4\pi }\iint \log |x-y|\,d\mu (x)\,d\mu (y)+\frac{1}{4\pi }%
\int |x|^{2}\,d\mu (x),  \label{eq:Einf-numerics}
\end{equation}%
the repulsion carrying $\tfrac{1}{2}G=-\tfrac{1}{4\pi }\log $ and the
confinement $\omega _{\ast }\!\int V=\tfrac{1}{2\pi }\!\int \tfrac{1}{2}%
|x|^{2}=\tfrac{1}{4\pi }\!\int |x|^{2}$. We recall that as shown in Section~\ref{sec:circle-example}, the single centered circle
minimizing~\eqref{eq:Einf-numerics} among centered circles sits at $r_{\ast }=1/\sqrt{2}$ with
energy $E(r_{\ast })=\tfrac{1}{8\pi }(1+\log 2)\approx 0.06737$.

\subsection{Two concentric circles}

\label{subsec:two-circles}

With $\nu _{r}$ as in Section~\ref{sec:circle-example}, consider the family of two-circle measures
\begin{equation*}
\mu =m_{1}\,\nu _{r_{1}}+m_{2}\,\nu _{r_{2}},\qquad m_{1}+m_{2}=1,\qquad
0<r_{1}<r_{2}.
\end{equation*}%
The only building block needed is the logarithmic interaction of two
concentric uniform circles. Writing $|r_{a}e^{i\alpha }-r_{b}|=\max (r_{a},r_{b})\,|e^{i\alpha }-q|$, after the reflection $\alpha \mapsto -\alpha $ when $r_{a}<r_{b}$, with $q=\min (r_{a},r_{b})/\max (r_{a},r_{b})\leq 1$, and using $\tfrac{1}{2\pi }\int_{0}^{2\pi }\log |e^{i\alpha }-q|\,d\alpha =0$ for $q\leq 1$, one obtains the exact
value
\begin{equation}
I(r_{a},r_{b}):=\iint \log |x-y|\,d\nu _{r_{a}}(x)\,d\nu _{r_{b}}(y)=\log
\max (r_{a},r_{b}).  \label{eq:two-circle-block}
\end{equation}%
In particular the self-energies are $I(r_{1},r_{1})=\log r_{1}$ and $%
I(r_{2},r_{2})=\log r_{2}$, while the cross term (in either order, by
symmetry) is $I(r_{1},r_{2})=\log r_{2}$. Expanding the double integral over
$\mu $ and adding the trap $\int |x|^{2}\,d\mu
=m_{1}r_{1}^{2}+m_{2}r_{2}^{2} $, the energy~\eqref{eq:Einf-numerics}
becomes
\begin{equation}
\begin{split}
E(m_{1},m_{2},r_{1},r_{2})& =-\frac{1}{4\pi }\Bigl[m_{1}^{2}\log
r_{1}+m_{2}^{2}\log r_{2}+2m_{1}m_{2}\log r_{2}\Bigr] \\
& \quad +\frac{1}{4\pi }\bigl(m_{1}r_{1}^{2}+m_{2}r_{2}^{2}\bigr).
\end{split}
\label{eq:Einf-two-circle}
\end{equation}%
The cross term $2m_{1}m_{2}\log r_{2}$ is the mutual energy of the two
circles and depends only on the outer radius, a direct consequence of~%
\eqref{eq:two-circle-block}.

With $m_{2}=1-m_{1}$ to enforce a unit total mass, the three stationarity conditions for~%
\eqref{eq:Einf-two-circle} are
\begin{align*}
\partial _{r_{1}}E& =-\frac{1}{4\pi }\frac{m_{1}^{2}}{r_{1}}+\frac{1}{4\pi }%
2m_{1}r_{1}=0\;\Longrightarrow \;r_{1}^{2}=\frac{m_{1}}{2}, \\
\partial _{r_{2}}E& =-\frac{1}{4\pi }\frac{m_{2}^{2}+2m_{1}m_{2}}{r_{2}}+%
\frac{1}{4\pi }2m_{2}r_{2}=0\;\Longrightarrow \;r_{2}^{2}=\frac{m_{2}+2m_{1}%
}{2}=\frac{1+m_{1}}{2}, \\
\partial _{m_{1}}E& =-\frac{1}{4\pi }2m_{1}\log \frac{r_{1}}{r_{2}}+\frac{1}{%
4\pi }\bigl(r_{1}^{2}-r_{2}^{2}\bigr)=0\;\Longrightarrow \;2m_{1}\log \frac{%
r_{2}}{r_{1}}=r_{2}^{2}-r_{1}^{2}=\tfrac{1}{2}.
\end{align*}%
Substituting the first two relations into the third collapses the system to
a single transcendental equation in the inner mass,
\begin{equation}
m_{1}\,\log \frac{1+m_{1}}{m_{1}}=\frac{1}{2},  \label{eq:transcendental}
\end{equation}%
which has a unique root $m_{1}\in (0,1)$. Seeking solutions numerically, we find
\begin{equation*}
\begin{aligned}
m_{1}& \approx 0.397953,& m_{2}&\approx 0.602047, \\
r_{1}& =\sqrt{\tfrac{m_{1}}{2}}\approx 0.446068,&\text{and}\quad r_{2}&=\sqrt{\tfrac{1+m_{1}}{2}}\approx 0.836048,
\end{aligned}
\end{equation*}%
with energy $E\approx 0.061956$. The $3\times 3$ Hessian in $(m_{1},r_{1},r_{2})$
is positive definite, so this is a local minimum within the two-circle
family, and it lies strictly below the single circle,
\begin{equation*}
E^{\text{two}}=0.06196<E^{\text{one}}(r_{\ast })=0.06737.
\end{equation*}%
Note that both masses are of order one and neither degenerates to zero.

\subsection{Configurations with many circles}

Adding circles can lower the energy further. The minimizer of~\eqref{eq:Einf-numerics} over all probability measures is the Frostman equilibrium measure, here the uniform measure on the unit disk~\cite{ST97,SandierSerfaty2015}, whose energy is $3/(16\pi )\approx 0.05968$. Measures carried by many concentric circles, with radii and masses chosen to mimic this uniform density, approach it as the number of circles grows, while no measure carried by finitely many curves reaches it. The two-circle energy $0.06196$ lies between the energy $0.06737$ of the single circle and the minimal energy $0.05968$. These observations suggest that the curve-supported critical points of the limit energy form a much richer family than the single circular solution, organized by the number of concentric circles and by the distribution of mass among them. We do not pursue the many-circle problem here.

On the discrete side, numerical experiments suggest that relative equilibria of $N$ point vortices exist close to the two-circle critical point of Section~\ref{subsec:two-circles}, with the vortices split as the masses $m_{1}$ and $m_{2}$. For configurations of concentric circles with a vortex at the center, the discrete counterparts are nested polygons that bifurcate from the regular $N$-gon with a central vortex, as shown locally by Meyer and Schmidt~\cite{MS88} and globally by Garc\'{\i}a-Azpeitia and Ize~\cite{GI11}. A vortex at the center, or equivalently a change of confining potential, can also stabilize the polygon. For identical vortices, the regular $N$-gon with a vortex at its center is a nondegenerate local minimizer of the energy, up to rotation, for $4\leq N\leq 8$, as follows from the Hessian computed in~\cite{GI12}, while the $N$-gon alone is stable only for $N\leq 6$ (Section~\ref{subsec:linearization-circle}). By Theorem~\ref{thm:energy-instability}, no such stabilization occurs on a curve.

We conclude with two natural open questions that arise from this discussion and its examples. They concern the equilibria of the discrete energy~\eqref{eq:En} of our model. The first is whether the critical points on several circles are limits, as $N\rightarrow \infty $, of these equilibria, just as the single circle is the limit of the regular $N$-gons (Section~\ref{sec:circle-example}), and as the numerical experiments suggest for the two-circle critical point. The second is whether stable equilibria can be found for larger and larger $N$, close to these critical points, as more and more circles are added, just as one vortex at the center extends the stability of the regular $N$-gon from $N=6$ up to $N=8$.\medskip

\noindent
\textbf{Open access.}
 For the purpose of open access, the authors have applied a Creative Commons Attribution (CC BY) licence to any Author Accepted Manuscript
version arising from this submission.\smallskip

\noindent
\textbf{AI usage statement.}
AI tools were used in support of the present work to aid with literature research, and in devising some simple numerical experiments to support the formulation of the results. The authors have personally checked, verified and written up all results reported.

\end{document}